\documentclass[10pt,a4paper]{article}	
\usepackage{amsfonts}			
\usepackage{amsmath}			
\usepackage{amsthm}			
\usepackage{amssymb}			
\usepackage{makeidx}
\usepackage{graphicx}			
\usepackage{enumerate}			
\usepackage{mathrsfs}			
\usepackage{hyperref}
\usepackage{pstricks}			
\usepackage[all,knot,arc]{xy}		
\usepackage{rotating}
\usepackage{array}

\graphicspath{{../immagini/}}

\newcommand{\executeiffilenewer}[3]{%
\ifnum\pdfs

trcmp{\pdffilemoddate{#1}}%
{\pdffilemoddate{#2}}>0%
{\immediate\write18{#3}}\fi%
}

\theoremstyle{plain}
\newtheorem{mythm}{Theorem}[section]
\newtheorem{mylem}[mythm]{Lemma}
\newtheorem{myprop}[mythm]{Proposition}

\theoremstyle{definition}
\newtheorem{mydef}[mythm]{Definition}
\newtheorem{myoss}[mythm]{Remark}
\newtheorem{myes}[mythm]{Example}

\newtheorem*{mystep}{}

\theoremstyle{plain}
\newtheorem{myintrothm}{Theorem}

\newtheorem{myintrocor}[myintrothm]{Corollary}

\theoremstyle{definition}
\newtheorem*{myintrooss}{Remark}

\newcommand{\mc}[1]{\mathcal{#1}}

\newcommand{\mg}[1]{\mathfrak{#1}}

\newcommand{\on}[1]{\operatorname{#1}}

\newcommand{\myemph}[1]{\emph{#1}}

\newcommand{\nN}{\mathbb{N}}

\newcommand{\nQ}{\mathbb{Q}}

\newcommand{\nC}{\mathbb{C}}
\newcommand{\nD}{\mathbb{D}}

\newcommand{\nK}{\mathbb{K}}

\newcommand{\nmat}[3]{{\mc{M}(#1 \times #2, #3)}}
\newcommand{\nind}[1]{{\mc{N}({#1})}}
\newcommand{\nindd}[2]{{\mc{N}_{#2}({#1})}}
\newcommand{\crit}[1]{{\mc{C}(#1)}}
\newcommand{\critinf}[1]{\crit{{#1}^\infty}}
\newcommand{\one}{{1\hspace*{-0.8 mm}\mathrm{l}}}

\newcommand{\x}{{w}}				
\newcommand{\ve}[1]{{#1}}			
\newcommand{\vx}{\x}				

\newcommand{\xu}{{u}}
\newcommand{\xv}{{v}}
\newcommand{\xx}{{x}}
\renewcommand{\xy}{{y}}
\newcommand{\xz}{{z}}
\newcommand{\xt}{{t}}

\newcommand{\vxu}{\ve{\xu}}
\newcommand{\vxv}{\ve{\xv}}
\newcommand{\vxx}{\ve{\xx}}
\newcommand{\vxy}{\ve{\xy}}
\newcommand{\vxz}{\ve{\xz}}
\newcommand{\vxt}{\ve{\xt}}

\newcommand{\X}{{x}}
\newcommand{\Y}{{y}}
\newcommand{\Z}{{z}}
\newcommand{\T}{{t}}

\newcommand{\vi}[1]{{#1}}		
\newcommand{\mi}[1]{{#1}}		
\newcommand{\vpp}[1]{{\uppercase{#1}}}	
\newcommand{\cm}[3]{{{#1}_{#2}^{#3}}}		
\newcommand{\cc}[3]{{#1}^{#2}_{#3}}		

\newcommand{\lhs}{\mathbb{I}}
\newcommand{\rhs}{\mathbb{II}}
\newcommand{\eq}{\mathbb{E}}
\newcommand{\cond}[1]{{\on{Cond}_{#1}}}

\newcommand{\wa}[2]{\on{weight}(#1,#2)}
\newcommand{\wb}[1]{{#1}}

\newcommand{\second}{{\prime\prime}}

\newcommand{\abs}[1]{\left|#1\right|}

\newcommand{\norm}[1]{\left\|#1\right\|}

\newcommand{\comp}{\circ}
\newcommand{\lot}{\on{l.o.t.}}
\newcommand{\lllot}[3]{\lot_{#1, #2}(#3)}
\newcommand{\llot}[2]{\lot_{#1}(#2)}

\newcommand{\kronecker}[2]{\delta_{#1}^{#2}}

\newcommand{\idmat}{{\on{Id}}}

\newcommand{\sumaa}[2]{\hspace{-#2}\sum_{#1}\hspace{-#2}}

\newenvironment{mydim}[0]		
{\begin{proof}}
{\end{proof}}

\newenvironment{mytext}[0]{}{}		

\newenvironment{mymatrix}[1]
{
\left(
\begin{array}{*{#1}{c}}
}
{
\end{array}
\right)
}

\newenvironment{mybmatrix}[1]
{
\left(
\begin{array}{*{#1}{c|}c}
}
{
\end{array}
\right)
}

\newcommand{\mydefin}[1]{\textit{#1}} 

\newcommand{\mysubsect}[1]
{
\subsection{#1}
}

\newcommand{\mysect}[1]
{
\section{#1}
\markright{\hfill \thesection\ #1 \hfill}
}

\newcommand{\mysecta}[1]
{
\section*{#1}
\markright{\hfill \thesection\ #1 \hfill}
}

\newcommand{\myacknowledgements}
{
\subsubsection*{Acknowledgements}
}

\newenvironment{myabstract}[1]
{\begin{center}\small{\textbf{#1}}

\vspace{0.25cm}
\begin{minipage}[l]{10.4cm}\begin{small}}
{\end{small}\end{minipage}\end{center}}

\title{Contracting rigid germs in higher dimensions.}
\author{Matteo Ruggiero}
\date{}

\makeindex

\begin{document}

\maketitle

\begin{myabstract}{Abstract}
Following Favre, we define a holomorphic germ $f:(\nC^d, 0) \rightarrow (\nC^d, 0)$ to be rigid if the union of the critical set of all iterates has simple normal crossing singularities. 
We give a partial classification of contracting rigid germs in arbitrary dimensions up to holomorphic conjugacy.
Interestingly enough, we find new resonance phenomena involving the differential of $f$ and its linear action on the fundamental group of the complement of the critical set.
\end{myabstract}

\pagestyle{plain}				


\mysecta{Introduction}
In this paper, we are concerned with the problem of analytic and formal classifications of contracting holomorphic germs at the origin in $\nC^d$, i.e., holomorphic germs $f:(\nC^d, 0) \rightarrow (\nC^d,0)$ such that every eigenvalue $\lambda$ of the differential $df_0$ at $0$ satisfies $0 \leq \abs{\lambda} < 1$.
The case of locally invertible maps is treated in detail in the literature (see, e.g., \cite{sternberg:localcontractionstheorempoincare}, \cite{rosay-rudin:holomorphicmaps} or \cite[Chapter $4$]{berteloot:methodeschangementechelles}).
Such a map is not necessarily linearizable, but is analytically conjugated to a polynomial normal form that involves only resonant monomials.
In particular, the analytic and formal classifications coincide. 
When the map is not invertible, the situation is far more complicated since the topological type of the critical set and its images are formal (hence analytic) invariants of conjugacy.

To get around this problem, a natural class of maps has been introduced in \cite{favre:rigidgerms} and was referred to as \myemph{rigid germs}.
A rigid germ $f:(\nC^d, 0) \rightarrow (\nC^d,0)$ is a holomorphic germ for which the \mydefin{generalized critical set} $\critinf{f}=\bigcup_{n\in \nN} f^{-n}\crit{f}$, where $\crit{f}$ denotes the critical set of $f$, has simple normal crossing singularities at the origin.

Favre gave a complete classification of contracting rigid germs in dimension $2$ (when $\critinf{f}$ is also totally invariant), and proved these germs were conjugated to a polynomial (or rational) normal form.
This classification has very interesting applications to the study of a special class of non-K\"ahler compact complex surfaces: Kato surfaces (see, e.g., \cite{dloussky:phdthesis}, \cite{dloussky:classgermesholocontr}, \cite{dloussky-oeljeklaus-toma:class70surfaces}).
The importance of this class was further emphasized by the work of \cite{favre-jonsson:eigenval} and \cite{ruggiero:rigidification}, since any holomorphic two-dimensional germ is birationally conjugated to a rigid germ.

In this article, we explore the classification of contracting rigid germs in higher dimensions.

In this setting, we shall exhibit new resonance phenomena involving the differential of $f$ at $0$ and its linear action on the fundamental group of the complement of the generalized critical set.
We shall then give some partial results on the classification of contracting rigid germs (see Theorem \ref{thm:rigidsecondaryresonance} and Theorem \ref{thm:rigidaffine}).

We shall also show (see Examples \ref{es:anycurve} and \ref{es:manyimages}) how the complexity of the geometry of the images of $\critinf{f}$ by $f$ and its iterates makes it impossible to get an explicit full classification.

A general motivation for studying contracting rigid germs in higher dimensions comes from the close relationship between these objects and special non-K\"ahler manifolds introduced by Kato (see, e.g., \cite{kato:cptcplxmanifoldsGSS}, \cite{kato:cpctcplx3fldlines}, \cite{oeljeklaus-renaud:3foldpolyauto}).
We shall return to this in a later work.

\medskip

The first natural invariant for holomorphic germs $f:(\nC^d, 0) \rightarrow (\nC^d, 0)$ is given by the differential $df_0$ at $0$.
In particular the number of non-zero eigenvalues of $df_0$, or equivalently the rank of $df_0^d$, is also invariant under iterations.

For contracting rigid germs one can consider a second natural invariant related to the (generalized) critical set.
Let $W^1, \ldots, W^q$ be the irreducible components of $\critinf{f}$: since $f$ is rigid, they are smooth and intersect transversely at $0$; moreover, since $\critinf{f}$ is backward invariant, for every $k=1, \ldots, q$ we have
$$
f^*W^k=\sum_{l=1}^q \cm{a}{l}{k} W^l
$$
for suitable $\cm{a}{l}{k} \in \nN$.
We define the \mydefin{internal action} of $f$ to be the matrix $A=A(f):=(\cm{a}{l}{k})$.
It can be understood geometrically since $A$ represents the action of $f$ on the fundamental group $\pi_1(\Delta^d \setminus \critinf{f})$, where $\Delta^d$ is a small open polydisc centered in $0$.

If $W^k=\{\vxu^k=0\}$ is a periodic component for $f^*$, i.e., ${(f^*)}^\eta W^k=W^k$ for a suitable $\eta \in \nN^*$, then $\left.\frac{\partial}{\partial \vxu^k}\right|_{0}$ defines an eigenvector for $df_0$ associated to a non-zero eigenvalue.
These eigenvalues are responsable for (part of) the classical resonances, given by the Poincar\'e-Dulac theorem.

For sake of simplicity, assume all periodic components are fixed and the nilpotent part of $df_0$ vanishes.
Observe that this can always be achieved replacing $f$ by a suitable iterate.

We shall also assume that $A$ is \myemph{injective}.
We observe that in dimension $2$, one can always semi-conjugate a rigid germ $f$ to another one $g$ satisfying this condition by \cite[Proposition 1.4]{favre:rigidgerms}, \cite[Theorem 5.1]{favre-jonsson:eigenval}, \cite[Remark 4.8]{ruggiero:rigidification} (see also Remark \ref{oss:noninjectiveinternalaction}): there exists a (not necessarily invertible) holomorphic germ $\Phi:(\nC^2,0) \rightarrow (\nC^2,0)$ such that $\det d\Phi_0 \not \equiv 0$ and $\Phi \circ f = g \circ \Phi$.

Classical resonances appear as algebraic relations between the eigenvalues of $df_0$.
The analysis of these resonances leads to the Poincar\'e-Dulac theorem.
Studying contracting rigid germs, a second kind of resonances appears, involving algebraic relations between (non-zero) eigenvalues of $df_0$ and eigenvalues of the internal action $A(f)$.

Let us denote by $\ve{\lambda} \in (\nD^*)^s$ the vector of non-eigenvalues of $df_0$ (where $\nD^*$ denotes the punctured unitary disc in $\nC$), and pick some coordinates $\vx=(\vxx, \cdot)$ such that $df_0=\on{Diag}(\ve{\lambda}, 0)$.

Set $\vxx=(\xx^1, \ldots, \xx^s)$ and let $\vi{n}=(n_1, \ldots, n_s)\in \nN^s$ be a multi-index with $\abs{n}:=n_1 + \ldots + n_s \geq 1$.
Then a monomial $\vxx^\vi{n}:=\prod_{k=1}^s(\xx^k)^{n_k}$ is \mydefin{secondary resonant}\label{pag:defintrosecondaryresonant} if and only if $\ve{\lambda}^\vi{n}$ is an eigenvalue for $\mi{A}$.

We can now state our main result. 

\begin{myintrothm}\label{thm:rigidresonantpart}
Let $f:(\nC^d, 0) \rightarrow (\nC^d,0)$ be a contracting rigid germ with injective internal action.
Suppose that all periodic components of $\critinf{f}$ are fixed, and the nilpotent part of $df_0$ vanishes.
Then $f$ is holomorphically conjugated to a map of the form
\begin{equation}\label{eqn:rigidresonantpart}
(\vxx,\vxy,\vxz) \mapsto \Big(\ve{\sigma}_{\on{PD}}(\vxx), \ve{\beta} \vxx^\mi{E}\vxy^\mi{D} \big(\one+\ve{g}(\vxx)\big), \ve{h}(\vxx,\vxy,\vxz) \Big)\mbox,
\end{equation}
where
\begin{itemize}
\item $\vxx \in \nC^s$, $\vxy \in \nC^p$, $\vxz \in \nC^{d-(s+p)}$, $\ve{\beta} \in (\nC^*)^p$ and $\one = (1, \ldots, 1)$;
\item $\mi{E} \in \nmat{s}{p}{\nN}$ and $\mi{D} \in \nmat{p}{p}{\nN}$ are matrices with $\det D \neq 0$;
\item $\ve{\sigma}_{\on{PD}}:(\nC^s,0)\rightarrow (\nC^s, 0)$ is a contracting invertible germ in Poincar\'e-Dulac normal form (hence a polynomial), and $df_0=d(\sigma_{\on{PD}})_0 \oplus 0$;
\item $\{\vxy^\one=0\} \subseteq \critinf{f} \subseteq \{\vxx^\one \vxy^\one = 0\}$;
\item $\ve{g}:(\nC^s, 0) \rightarrow (\nC^p,0)$ is a polynomial map that contains only secondary resonant monomials;
\item $\ve{h}:(\nC^d, 0) \rightarrow (\nC^{d-(s+p)},0)$ is a holomorphic map such that $dh_0=0$.
\end{itemize}
\end{myintrothm}
\begin{myintrooss}
Theorem \ref{thm:rigidresonantpart} still holds if we replace $\nC$ by a (possibly non-archimedean, not algebraically closed) complete metrized field $\nK$ of $\on{char}(\nK)=0$, provided that the eigenvalues of $df_0$ belong to $\nK$.
See Remarks \ref{oss:nonarchistart}, \ref{oss:nonarchilinear}, \ref{oss:nonarchijordan}, \ref{oss:nonarchiprimary}, \ref{oss:nonarchisecondary}, \ref{oss:nonarchiaffine} for further details.

Notice that Theorem \ref{thm:rigidresonantpart} does not hold over fields of positive characteristic, already for $d=p=1$ (see Remark \ref{oss:nonarchisecondary} for further details).
\end{myintrooss}

To read \eqref{eqn:rigidresonantpart}, we set $\vxx^\mi{E}=\Big(\big(\vxx^\mi{E}\big)^1, \ldots, \big(\vxx^\mi{E}\big)^p\Big) \in \nC^p$, with
$$
\big(\vxx^\mi{E}\big)^k=\prod_{l=1}^s (\xx^l)^{\cm{e}{l}{k}}\mbox,
$$
where $\vxx=(\xx^1, \ldots, \xx^s)$ and $\mi{E}=(\cm{e}{l}{k})$, and analogously for other similar expressions.
Moreover, if $\ve{\alpha}=(\alpha^1, \ldots, \alpha^p) \in \nC^p$ and $\vxy=(\xy^1, \ldots, \xy^p) \in \nC^p$, we shall write
$$
\ve{\alpha}\vxy=(\alpha^1 \xy^1, \ldots, \alpha^p \xy^p) \in \nC^p\mbox.
$$

This normal form has several features.
The first part $\vxx \circ \ve{f}$ depends only on $\vxx$, and defines a polynomial biholomorphism $\sigma_{\on{PD}}:\nC^s \rightarrow \nC^s$.
The second part $\vxy \circ \ve{f}$ depends only (polynomially) on $\vxx$ and (monomially) on $\vxy$.
We also get the following corollary.
\begin{myintrocor}
Let $f:(\nC^d, 0) \rightarrow (\nC^d,0)$ be a contracting rigid germ satisfying the hypotheses of Theorem \ref{thm:rigidresonantpart}.
Then $f$ preserves (at least) $s+1$ smooth foliations $\mc{F}_1, \ldots, \mc{F}_s$, $\mc{G}$.
The foliation $\mc{F}_k$ has codimension $k$ for $k=1, \ldots, s$, while $\mc{G}$ has codimention $s+p$.
Moreover, $\mc{F}_l$ is a subfoliation of $\mc{F}_k$ for every $l \geq k$, and $\mc{G}$ is a subfoliation of $\mc{F}_s$.
\end{myintrocor}

In the following we shall deal with rigid germs without the assumptions on the fixed components and the nilpotent part.
We shall then prove Theorem \ref{thm:rigidsecondaryresonance} as a generalization of Theorem \ref{thm:rigidresonantpart}.
We shall also study the special case of a rigid germ $f:(\nC^d, 0) \rightarrow (\nC^d, 0)$ with $s+p=d-1$, where $s$ is the number of non-zero eigenvalues of $df_0$, and $p$ is the number of non-periodic irreducible components of $\critinf{f}$ (see Theorem \ref{thm:rigidaffine}).
In particular we get the classifications for rigid germs for which $\critinf{f}$ has either $(d-1)$ or $d$ irreducible components.

These results solve the classification of rigid germs in dimension $3$, but for the case $s+p=1$.
We shall show (see Examples \ref{es:anycurve} and \ref{es:manyimages}) how an explicit classification of rigid germs in this case is not possible.

\medskip

To prove Theorem \ref{thm:rigidresonantpart}, we first apply Poincar\'e-Dulac normalization techniques.
We then get a germ $f:(\nC^d,0)\rightarrow(\nC^d,0)$ whose first $s$ coordinates are given by a contracting invertible polynomial $\ve{\sigma}_{\on{PD}}:(\nC^s,0)\rightarrow (\nC^s, 0)$ in Poincar\'e-Dulac normal form.
We then use the rigid assumption to get \eqref{eqn:rigidresonantpart} with $\ve{g}=\ve{g}(\vxx, \vxy, \vxz)$ that a priori depends on all coordinates.
We finally conjugate again to get $\ve{g}=\ve{g}(\vxx)$ that depends only on $\vxx$.
The study of resonances in this case will allow us to get $\ve{g}$ polynomial.

Conjugations that maintain the Poincar\'e-Dulac normal form are called \myemph{renormalizations}.
See for example \cite{abate-tovena:formalnormalformsholotangid} for a renormalization process in the tangent-to-the-identity case, \cite{abate-raissy:renormalization} for a general procedure for formal renormalizations, or \cite{raissy:torusactionnormprob} for other techniques to study the convergence of the Poincar\'e-Dulac normalization.

At every step, we first deal with the formal conjugacy.
To expand in formal power series compositions of maps, we shall need to introduce matrices of indices.
This will allow to fully understand the combinatorial structure of the formal problem.
Dealing with the general case (where the components of $\critinf{f}$ are not necessarily fixed, and $df_0$ has a nilpotent part) will make the combinatorial structure even more intricate.

Once the formal normal form is achieved, we solve the convergence problem as in the Poincar\'e-Dulac result.

\medskip

This paper is organized as follows.
In Section $1$ we prove some preparatory lemmas.
In Section $2$ we prove a generalization of the Poincar\'e-Dulac Theorem suited for contracting rigid germs.
In Section $3$ we define secondary resonances and prove Theorem \ref{thm:rigidresonantpart} in the general case.
In Section $4$ we deal with the special case $s+p=d-1$.
In Section $5$ we specify all results to dimension $3$.

\myacknowledgements
I would like to thank Marco Abate for many suggestions, and Charles Favre for very fruitful discussions.
This work was partially carried out during my PhD at Scuola Normale Superiore in Pisa, and visiting the \'Ecole Polytechnique in Paris. I would like to thank both institutions for their hospitality.
The grants of the project ANR Berko during my visiting period in Paris were fundamental and greatly acknowledged.

\mysect{Linear Part}
Let $f: (\nC^d, 0) \rightarrow (\nC^d, 0)$ be a contracting rigid germ, and $W^1, \ldots, W^q$ be the irreducible components of $\critinf{f}$.
It is natural to choose coordinates $\vx=(\x^1, \ldots, \x^d)$ such that $W^k=\{\x^k = 0\}$ for $k=1, \ldots, q$.
Moreover, as anticipated in the introduction, we want to split irreducible components between periodic and non-periodic ones with respect to the action of $f^*$.
Up to permuting coordinates, we can then suppose that the matrix $A=A(f)$ is of the form
\begin{equation}\label{eqn:definternalaction}
\mi{A}=
\begin{mymatrix}{2}
\mi{B} & \mi{C} \\
\mi{0} & \mi{D}
\end{mymatrix}
\mbox,
\end{equation}
where $\mi{B} \in \nmat{r}{r}{\nN}$ for a suitable $0 \leq r \leq q$ is a permutation matrix, $\mi{C} \in \nmat{r}{p}{\nN}$ with $p=q-r$, and $\mi{D} \in \nmat{p}{p}{\nN}$.

Since $\critinf{f}$ is backward $f$-invariant, in the chosen coordinates we can write
\begin{equation}\label{eqn:rigidstart}
(\vxu, \vxy, \vxt) \mapsto \Big(\ve{\alpha} \vxu^{\mi{B}} \big(\one + \ve{\theta}(\vxu, \vxy, \vxt)\big), \ve{\beta} \vxu^\mi{C}\vxy^\mi{D} \big(\one+\ve{g}(\vxu,\vxy,\vxt)\big), \ve{k}(\vxu,\vxy,\vxt) \Big)\mbox,
\end{equation}
with
\begin{itemize}
\item $\vxu \in \nC^r$, $\vxy \in \nC^p$ and $\vxt \in \nC^{d-q}$;
\item $\ve{\alpha} \in (\nC^*)^r$ and $\ve{\beta} \in (\nC^*)^p$;
\item $\ve{\theta}:(\nC^d,0) \rightarrow (\nC^r, 0)$, $\ve{g}:(\nC^d,0) \rightarrow (\nC^p, 0)$, $\ve{k}:(\nC^d,0) \rightarrow (\nC^{d-q}, 0)$;
\item $\critinf{f}=\{\vxu^\one \vxy^\one=0\}$.
\end{itemize}

\begin{myoss}\label{oss:nonarchistart}
Observe that this reduction to maps of the form \eqref{eqn:rigidstart} is also valid over an arbitrary field.
\end{myoss}

\begin{myoss}
Suppose $q = d$.
If $\det{\mi{D}} \neq 0$, the condition $\critinf{f}=\{\vxu^\one \vxy^\one=0\}$ implies that the matrix $\mi{C}=(\cm{c}{i}{j})_{i,j}$ satisfies $\cm{c}{i}{1} + \cdots + \cm{c}{i}{p} \geq 1$ for every $i=1, \ldots, r$.
It is easy to check that every germ of the form \eqref{eqn:rigidstart} with $\mi{C}$ as above is a rigid germ.

If $\det{\mi{D}} = 0$, then every rigid germ can be written in the form \eqref{eqn:rigidstart}, but not every germ of the form \eqref{eqn:rigidstart} is a rigid germ.
For example, let us consider $f:(\nC^2, 0) \rightarrow (\nC^2, 0)$, given by 
$$
(\xy^1, \xy^2) \mapsto \big(\xy^1\xy^2 (1+\xy^1),\xy^1\xy^2 (1+\xy^2)\big)\mbox.
$$
Here
\begin{equation*}
\mi{A}=\mi{D}=
\begin{mymatrix}{2}
1 & 1 \\
1 & 1
\end{mymatrix}
\mbox,
\end{equation*}
hence $\det{\mi{D}}=0$, while $\det df = \xy^1 \xy^2 (\xy^1 + \xy^2 + 3 \xy^1\xy^2)$, hence $f$ is not rigid.

When $q < d$, a germ of the form \eqref{eqn:rigidstart} needs to satisfy suitable additional conditions (depending on $\mi{C}, \ve{k}$ and $\ve{g}$ if $\det{\mi{D}}=0$) to be rigid.
\end{myoss}

Starting from a germ of the form \eqref{eqn:rigidstart}, we would like to kill $\theta$ (up to holomorphic conjugacy). This is exactly the result of this section.

\begin{mythm}\label{thm:rigidlinearpart}
Let $\ve{f}:(\nC^d, 0) \rightarrow (\nC^d,0)$ be a contracting rigid germ.
Then $f$ is holomorphically conjugated to
\begin{equation}\label{eqn:rigidintlinearpart}
(\vxu, \vxy, \vxt) \mapsto \Big(\ve{\alpha} \vxu^{\mi{B}}, \ve{\beta} \vxu^\mi{C}\vxy^\mi{D} \big(\one+\ve{g}(\vxu,\vxy,\vxt)\big), \ve{k}(\vxu,\vxy,\vxt) \Big)\mbox,
\end{equation}
where
\begin{itemize}
\item $\vxu \in \nC^r$, $\vxy \in \nC^p$ and $\vxt \in \nC^{d-q}$;
\item $\ve{\alpha} \in (\nC^*)^r$ and $\ve{\beta} \in (\nC^*)^p$;
\item $\mi{B} \in \nmat{r}{r}{\nN}$ is a permutation matrix, $\mi{C} \in \nmat{r}{p}{\nN}$ and $\mi{D} \in \nmat{p}{p}{\nN}$;
\item $\ve{g}:(\nC^d,0) \rightarrow (\nC^p, 0)$, $\ve{k}:(\nC^d,0) \rightarrow (\nC^{d-q}, 0)$;
\item $\critinf{f}=\{\vxu^\one \vxy^\one=0\}$.
\end{itemize}
\end{mythm}

\begin{mytext}
Before proving this theorem, we need an easy notation Lemma and a Proposition.
\end{mytext}

\begin{mylem}\label{lem:logf^D}
Let $\ve{f}=(f_1, \ldots, f_r):(\nC^d, 0) \rightarrow (\nC^r,0)$ be an $r$-uple of formal power series, and let $\mi{D} \in \nmat{r}{s}{\nQ}$.
Then we have
$$
\log \big(\ve{f}^{\mi{D}}\big)=(\log \ve{f}) \mi{D}\mbox,
$$
where $\log$ here means that we are taking the $\log$ coordinate by coordinate.
\end{mylem}
\begin{mydim}
It easily follows from a direct computation.
\end{mydim}

\begin{myprop}\label{prop:takingthelog}
Let $(\ve{f}_n)_n$ be a sequence of $r$-uples of formal power series, and let $(\mi{D}_n)_n$ be a sequence of matrices in $\nmat{r}{s}{\nQ}$ (with $s \geq 1$). Then
$$
\prod_n \big(\one+\ve{f}_n\big)^{\mi{D}_n}
$$
converges if and only if
$$
\sum_n \ve{f}_n \mi{D}_n 
$$
does.
\end{myprop}
\begin{mydim}
It follows from Lemma \ref{lem:logf^D} and the analogous result in dimension one, taking the $\log$ of the absolute value.
\end{mydim}

\begin{mydim}[Proof of Theorem \ref{thm:rigidlinearpart}]
We can suppose that $\ve{f}$ is of the form \eqref{eqn:rigidstart}.

We would like to find a conjugacy between $\ve{f}$ and a germ $\ve{\widetilde{f}}:(\nC^d, 0) \rightarrow (\nC^d, 0)$ of the form \eqref{eqn:rigidintlinearpart} (with $\ve{\widetilde{g}}, \ve{\widetilde{k}}$ replaced by some holomorphic maps $\ve{g}$ and $\ve{k}$ respectively).

Let us consider a local diffeomorphism of the form
$$
\ve{\Phi}(\vxu, \vxy, \vxt)=\Big(\vxu \big(\one+ \ve{\phi}(\vxu, \vxy, \vxt)\big), \vxy,\vxt \Big)\mbox,
$$
where $\ve{\phi}:(\nC^d, 0) \rightarrow (\nC^r,0)$.
Set 
$$
\ve{\Phi}_N(\vxu, \vxy, \vxt)=\Big(\vxu \big(\one + \ve{\phi}_N(\vxu, \vxy, \vxt)\big), \vxy, \vxt \Big)\mbox,
$$
with
$$
\one+\ve{\phi}_N(\vx)=\prod_{n=1}^N \big(\one+\ve{\theta} \circ \ve{f}^{\comp n-1}(\vx)\big)\mbox,
$$
where $\vx=(\vxu, \vxy, \vxt)$.
We have
$$
\vxu \circ \ve{\Phi}_N \circ \ve{f} = \vxu \circ \ve{\widetilde{f}} \circ \ve{\Phi}_{N+1}\mbox.
$$
Indeed
\begin{align*}
\vxu \circ \ve{\Phi}_N \circ \ve{f} (\vx) &= \ve{\alpha} \vxu \big(\one + \ve{\theta}(\vx)\big) \prod_{n=1}^N \big(\one+\ve{\theta} \circ \ve{f}^{\comp n}(\vx)\big)\mbox,\\
\vxu \circ \ve{\widetilde{f}} \circ \ve{\Phi}_{N+1} (\vx) &= \ve{\alpha} \vxu \prod_{n=1}^{N+1} \big(\one+\ve{\theta} \circ \ve{f}^{\comp n-1}(\vx)\big)\mbox,
\end{align*}
which are equivalent expressions.

Let us prove that $\ve{\Phi}_N$ converges to a holomorphic germ $\ve{\Phi}=\ve{\Phi}_\infty$.
Thanks to Proposition \ref{prop:takingthelog}, we just have to prove that
$$
\sum_{n=0}^\infty \ve{\theta} \circ \ve{f}^{\comp n}
$$
converges in a neighborhood of $0$.
Since $\ve{\theta}(\ve{0})=\ve{0}$, we have that there exists $M > 0$ such that $\norm{\ve{\theta}(\vx)} \leq M \norm{\vx}$,
while being $\ve{f}$ contracting, there exists $0 < \Lambda < 1$ such that $\norm{\ve{f}^{\comp n}(\vx)} \leq \Lambda^n \norm{\vx}$, both estimates for $\norm{\vx}$ small enough.
Then we have
$$
\norm{\sum_{n=0}^\infty \big(\ve{\theta} \circ \ve{f}^{\comp n}(\vx)\big)} \leq \sum_{n=0}^\infty \norm{\ve{\theta} \circ \ve{f}^{\comp n}(\vx)} \leq \sum_{n=0}^\infty M \Lambda^n \norm{\vx} = \frac{M}{1-\Lambda} \norm{\vx} < +\infty\mbox,
$$
and hence $\ve{\Phi}:(\nC^d,0)\rightarrow (\nC^d, 0)$ is a holomorphic invertible map, that satisfies the conjugacy relation
\begin{equation}\label{eqn:linearpartconjrel}
\ve{\Phi} \circ \ve{f}=\ve{\widetilde{f}} \circ \ve{\Phi}
\end{equation}
for the first coordinate $\vxu$.

We can just define
\begin{align*}
\one+\ve{\widetilde{g}}(\vx)&=\left(\frac{\one + \ve{g}}{(\one+\ve{\phi})^C}\right) \circ \ve{\Phi}^{-1}(\vx)\mbox,\\
\ve{\widetilde{k}}(\vx)&=\ve{k}\circ \ve{\Phi}^{-1}(\vx)\mbox,
\end{align*}
to have \eqref{eqn:linearpartconjrel} satisfied for all coordinates.
\end{mydim}
  
\begin{myoss}\label{oss:nonarchilinear}
Observe the arguments of the proof of Theorem \ref{thm:rigidlinearpart} are also valid over any complete metrized field of characteristic $0$.
\end{myoss}

\mysect{Primary Resonances}
\mysubsect{Resonance Relation}

\begin{mytext}
Considering the differential $df_0$ at $0$ for a map $f:(\nC^d, 0) \rightarrow (\nC^d, 0)$ of the form \eqref{eqn:rigidintlinearpart}, we get
$$
df_0=
\begin{mymatrix}{3}
\on{Diag}(\ve{\alpha}) \mi{B}^T & \mi{0} & \mi{0} \\
\ast & \mi{0} & \mi{0}\\
\ast & \ast & \Big. \left.\frac{\partial \ve{k}}{\partial \vxt}\right|_0 \Big.
\end{mymatrix}
\mbox,
$$
where $\on{Diag}(\ve{\alpha})$ is the diagonal matrix of entries $\ve{\alpha}=(\alpha^1, \ldots, \alpha^r)$, and $^T$ denotes the transposition.

Some of the non-zero eigenvalues of $df_0$ arise from the block $\on{Diag}(\ve{\alpha}) \mi{B}^T$, but they can arise also from $\left.\frac{\partial \ve{k}}{\partial \vxt}\right|_0$.
We can change coordinates to have $\left.\frac{\partial \ve{k}}{\partial \vxt}\right|_0$ in Jordan normal form, and split the coordinate $\vxt$ in $(\vxv, \vxz)$, with $\vxv \in \nC^e$ and $\vxz \in \nC^{d-(q+e)}$, to have the diagonal part equal to $(\ve{\mu}, \ve{0})$, where $\mu \in (\nD^*)^e$ is the vector of non-zero eigenvalues of $df_0$ that do not arise from $\on{Diag}(\ve{\alpha}) \mi{B}^T$.
With these new coordinates we can write $f$ as follows:
\begin{equation}\label{eqn:rigidsplitlinearpart}
(\vxu, \vxv, \vxy, \vxz) \mapsto \Big(\ve{\alpha} \vxu^{\mi{B}}, \ve{\mu} \vxv + \ve{\rho}(\vxu,\vxv,\vxy,\vxz), \ve{\beta} \vxu^\mi{C}\vxy^\mi{D} \big(\one +  \ve{g}(\vxu,\vxv,\vxy,\vxz)\big), \ve{h}(\vxu,\vxv,\vxy,\vxz) \Big)\mbox,
\end{equation}
where
\begin{itemize}
\item $\vxu \in \nC^r$, $\vxv \in \nC^e$ with $e=s-r$, $\vxy \in \nC^p$ and $\vxz \in \nC^{d-(s+p)}$;
\item $\ve{\alpha} \in (\nC^*)^r$, $\ve{\mu} \in (\nD^*)^e$ and $\ve{\beta} \in (\nC^*)^p$;
\item $\mi{B} \in \nmat{r}{r}{\nN}$ is a permutation matrix, $\mi{C} \in \nmat{r}{p}{\nN}$ and $\mi{D} \in \nmat{p}{p}{\nN}$;
\item $\ve{\rho}:(\nC^d,0) \rightarrow (\nC^e, 0)$, $\ve{g}:(\nC^d,0) \rightarrow (\nC^p, 0)$ and $\ve{h}:(\nC^d,0) \rightarrow (\nC^{d-(s+p)}, 0)$;
\item $\critinf{f}=\{\vxu^\one \vxy^\one=0\}$;
\item $\mu \subseteq \on{Spec}(df_0) \setminus \{0\}$;
\item $\ve{\rho}|_{\{\vxu=\vxy=\vxz=\ve{0}\}}$ and $\ve{h}|_{\{\vxu=\vxv=\vxy=\ve{0}\}}$ have nilpotent linear part.
\end{itemize}

\begin{myoss}\label{oss:nonarchijordan}
Over an arbitrary field $\nK$, this argument works as soon as all the eigenvalues of $df_0$ belong to $\nK$.
In particular, it always works if $\nK$ is algebraically closed.
\end{myoss}

We want now to kill as many coefficients of $\ve{\rho}$ (expanded in formal power series) as we can. As in the case of attracting invertible germs, some formal obstructions appear. We shall call them \myemph{primary} resonances to make a distinction with \myemph{secondary} resonances, that will be introduced in Definition \ref{def:secondaryresonances} (see also the introduction at page \pageref{pag:defintrosecondaryresonant}).
\end{mytext}

\begin{mydef}\label{def:primaryresonances}
Let $f:(\nC^d,0) \rightarrow (\nC^d,0)$ be a contracting rigid germ as in \eqref{eqn:rigidsplitlinearpart}, and let $\eta \in \nN^*$ be the order of $\mi{B}$.
A monomial $\vxu^{\vi{n}_\vxu} \vxv^{\vi{n}_\vxv}$ is called \mydefin{primary resonant} with respect to the $k$-th coordinate of $\vxv$ if it satisfies the Poincar\'e-Dulac resonance relation for $f^{\comp \eta}$, i.e., if
\begin{equation}\label{eqn:defprimaryresonance}
\ve{\xi}^{\vi{n}_\vxu}\ve{\mu}^{\eta \vi{n}_\vxv} = (\mu^k)^\eta\mbox,
\end{equation}
where $\ve{\xi} \in (\nD^*)^r$ is the vector of eigenvalues of $\left(\vxu \mapsto \ve{\alpha} \vxu^{\mi{B}}\right)^{\comp \eta}$ (counted with multiplicities), and $\mu=(\mu^1, \ldots, \mu^e)$.
\end{mydef}

\begin{myoss}
If $\eta=1$ in Definition \ref{def:primaryresonances} the resonance relation \eqref{eqn:defprimaryresonance} becomes
$$
\ve{\lambda}^\vi{n} = \lambda^{r+k}\mbox,
$$
where $\ve{\lambda}=(\lambda^1, \ldots, \lambda^s):=(\ve{\alpha}, \ve{\mu}) \in (\nD^*)^s$ is the vector of non-zero eigenvalues for $df_0$, and $\vi{n}=(\vi{n}_\vxu, \vi{n}_\vxv) \in \nN^s$.
\end{myoss}

\begin{myoss}
Let $f:(\nC^d,0) \rightarrow (\nC^d,0)$ be a contracting rigid germ as in \eqref{eqn:rigidsplitlinearpart}. Let us suppose for example that $W^k=\{\xu^k=0\}$ for $k=1, \ldots, \chi$ form a cycle of order $\chi$; i.e., the first $\chi$ coordinates of $f$ are of the form
$$
(\xu^1, \ldots, \xu^\chi)\mapsto(\alpha^1 \xu^2, \ldots, \alpha^{\chi-1} \xu^\chi, \alpha^\chi \xu^1)\mbox. 
$$
Taking the $\chi$-th iterate, we get
$$
(\xu^1, \ldots, \xu^\chi)\mapsto(\xi \xu^1, \ldots, \xi \xu^\chi)\mbox,\qquad \mbox{with } \xi = \prod_{k=1}^\chi \alpha^k\mbox. 
$$
In particular all $\left. \frac{\partial}{\partial \xu^k} \right|_0$ belong to the eigenspace of eigenvalue $\xi$ for $df_0^\chi$, and $\xi$ will have multiplicity (at least) $\chi$.
\end{myoss}

\begin{mytext}
The following lemma is a classical result for primary resonances in contracting germs (see, e.g., \cite[p. 467]{berteloot:methodeschangementechelles}).
\end{mytext}

\begin{mylem}
Let $\ve{f}:(\nC^d, 0) \rightarrow (\nC^d,0)$ be a contracting rigid germ written as in \eqref{eqn:rigidsplitlinearpart}.
Then there are only finitely-many primary resonant monomials.
\end{mylem}

\begin{myoss}
We notice that periodic non-fixed irreducible components of the generalized critical set of a contracting rigid germ $f:(\nC^d, 0) \rightarrow (\nC^d, 0)$ can appear only for $d \geq 3$, and primary resonances for $\mi{B} \neq \on{Id}$ can appear only for $d \geq 4$.
\end{myoss}

\mysubsect{Main Theorem}

\begin{mytext}
Our next goal is to kill all coefficients of $\ve{\rho}$ in \eqref{eqn:rigidsplitlinearpart} except for primary resonant monomials.
\end{mytext}

\begin{mythm}\label{thm:rigidprimaryresonance}
Let $\ve{f}:(\nC^d, 0) \rightarrow (\nC^d,0)$ be a contracting rigid germ.
Then $f$ is \myemph{analytically} conjugated to
\begin{equation}\label{eqn:rigidsplitprimaryresonance}
(\vxu, \vxv, \vxy, \vxz) \mapsto \Big(\ve{\alpha} \vxu^{\mi{B}}, \ve{\mu} \vxv + \ve{\rho}(\vxu,\vxv), \ve{\beta} \vxu^\mi{C}\vxy^\mi{D} \big(\one +  \ve{g}(\vxu,\vxv,\vxy,\vxz)\big), \ve{h}(\vxu,\vxv,\vxy,\vxz) \Big)\mbox,
\end{equation}
where
\begin{itemize}
\item $\vxu \in \nC^r$, $\vxv \in \nC^e$, $\vxy \in \nC^p$ and $\vxz \in \nC^{d-(s+p)}$;
\item $\ve{\alpha} \in (\nC^*)^r$, $\ve{\mu} \in (\nD^*)^e$ and $\ve{\beta} \in (\nC^*)^p$;
\item $\mi{B} \in \nmat{r}{r}{\nN}$ is a permutation matrix, $\mi{C} \in \nmat{r}{p}{\nN}$ and $\mi{D} \in \nmat{p}{p}{\nN}$;
\item $\ve{\rho}:(\nC^s,0) \rightarrow (\nC^e, 0)$, $\ve{g}:(\nC^d,0) \rightarrow (\nC^p, 0)$ and $\ve{h}:(\nC^d,0) \rightarrow (\nC^{d-(s+p)}, 0)$;
\item $\critinf{f}=\{\vxu^\one \vxy^\one=0\}$;
\item $\mu \subseteq \on{Spec}(df_0) \setminus \{0\}$ and $\ve{h}|_{\{\vxu=\vxv=\vxy=\ve{0}\}}$ has nilpotent linear part;
\item $\ve{\rho}$ is a polynomial map with only primary resonant monomials.
\end{itemize}
\end{mythm}
\begin{myoss}\label{oss:orderonmu}
For a contracting rigid germ $\ve{f}:(\nC^d, 0) \rightarrow (\nC^d,0)$ written as in \eqref{eqn:rigidsplitlinearpart}, up to permuting coordinates in $\vxv=(\xv^1, \ldots, \xv^e)$, we can order $\mu^1, \ldots, \mu^e$ such that
\begin{equation}\label{eqn:orderonmu}
1 > \abs{\mu^1} \geq \ldots \geq \abs{\mu^e} > 0\mbox.
\end{equation}
In this case a primary resonant monomial for the $k$-th coordinate is either of the form:
\begin{itemize}
\item $\vxu^{\vi{n}_\vxu} \vxv^{\vi{n}_\vxv}$ with $\vi{n}_\vxv=(n_{\vxv^1}, \ldots, n_{\vxv^e})$ such that $n_{\vxv^l}=0$ for $l \geq k$,
\item or $\xv^l$ for a suitable $1 \leq l \leq e$.
\end{itemize}
We shall first take care of the linear part, and then of higher order terms, to get $\ve{\rho}=(\rho^1, \ldots, \rho^e)$ (strictly) triangular, meaning precisely that $\rho^{k}$ depends only on $\vxu$ and $\xv^1, \ldots, \xv^{k-1}$ for every $k=1, \ldots, e$.
\end{myoss}

\begin{mydim}
We first prove in Step $1$ the formal counterpart of this theorem, and then we will deal with the convergence of the formal power series involved in Step $2$.

\begin{mystep}[Step $1$]

First, we can suppose that $f$ is of the form \eqref{eqn:rigidsplitlinearpart}.
Moreover we can suppose that $(\ve{\rho}, \ve{h})|_{\vxu=\vxy=0}$ has lower triangular linear part, and that $\mu$ satisfies \eqref{eqn:orderonmu}.

Then, up to linear conjugacy, we can suppose that the linear part of $\rho$ has only resonant monomials.
Indeed, we can consider a linear map of the form $\hat{L}:(\vxu,\vxv,\vxy,\vxz)\mapsto(L\vxu,\vxv,\vxy,\vxz)$ that conjugates $f$ with $\hat{f}= \hat{L} \circ f \circ \hat{L}^{-1}$ such that $\vxu \circ \hat{f}(\vxu, \vxv, \vxy, \vxz)=\hat{\alpha} \vxu$, where $\hat{\alpha} \in (\nD^*)^r$ is a vector of non-zero eigenvalues for $df_0$.
Then there is a linear change of coordinates $\hat{M}$ that conjugates $\hat{f}$ with a map whose linear part is in Jordan normal form, and $\hat{L}^{-1} \circ \hat{M} \circ \hat{L}$ is the wanted linear conjugacy.

Now we want to conjugate $f$ with a map $\widetilde{f}:(\nC^d, 0)\rightarrow (\nC^d, 0)$ of the form \eqref{eqn:rigidsplitprimaryresonance} (with $\widetilde{g}$, $\widetilde{\rho}$ and $\widetilde{h}$ instead of $g$, $\rho$ and $h$ respectively).

Set $\vx=(\vxu, \vxv, \vxy, \vxz)$, and consider a local diffeomorphism $\Phi:(\nC^d, 0) \rightarrow (\nC^d, 0)$ of the form
$$
\ve{\Phi}(\vx)=\big(\vxu, \ve{\phi}(\vx), \vxy, \vxz\big)\mbox,
$$
with $\phi:(\nC^d, 0)\rightarrow (\nC^e,0)$ a formal map such that $d\Phi_0=\idmat_d$ is tangent to the identity. 

Considering the conjugacy relation $\Phi \circ f = \widetilde{f} \circ \Phi$ for the coordinate $\vxv$, we have to solve
\begin{equation}\label{eqn:primaryeqntosolve}
\ve{\phi} \circ \ve{f}(\vx)
= \vxv \circ \ve{\Phi} \circ \ve{f} (\vx)
= \vxv \circ \ve{\widetilde{f}} \circ \ve{\Phi} (\vx)
= (\ve{\mu} \vxv + \ve{\widetilde{\rho}})\big(\vxu, \ve{\phi}(\vx)\big)
= \ve{\mu} \ve{\phi}(\vx) + \ve{\widetilde{\rho}}\big(\vxu, \ve{\phi}(\vx)\big)
\end{equation}
for suitable $\phi$ and $\widetilde{\rho}$.
Set
$$
\lhs(\vx):=\ve{\phi} \circ \ve{f}(\vx), \qquad
\rhs(\vx):= \ve{\mu} \ve{\phi}(\vx) + \ve{\widetilde{\rho}}\big(\vxu, \ve{\phi}(\vx)\big).
$$
We now expand in formal power series \eqref{eqn:primaryeqntosolve}, and solve it by defining recursively the coefficients of $\phi$ and $\widetilde{\rho}$.
Set 
\begin{itemize}
\item $\vxv=(\xv^1, \ldots, \xv^e)$;
\item $\ve{\rho}=(\rho^1, \ldots, \rho^e)$ and $\mu^k \xv^k + \rho^k(\vx)=\sum_{\vi{n}} \cc{\rho}{k}{\vi{n}}\vx^\vi{n}$ for $1 \leq k \leq e$;
\item $\ve{\widetilde{\rho}}=(\widetilde{\rho}^1, \ldots, \widetilde{\rho}^e)$ and $\mu^k \xv^k + \widetilde{\rho}^k(\vxu, \vxv)=\sum_{\vi{n}_\vxu, \vi{n}_\vxv} \cc{\widetilde{\rho}}{k}{(\vi{n}_\vxu,\vi{n}_\vxv)}\vxu^{\vi{n}_\vxu} \vxv^{\vi{n}_\vxv}$ for $1 \leq k \leq e$;
\item $\ve{\phi}=(\phi^1, \ldots, \phi^e)$ and $\phi^k(\vx)= \sum_{\vi{n}} \cc{\phi}{k}{\vi{n}}\vx^\vi{n}$ for $1 \leq k \leq e$;
\item $\lhs=\big(\lhs^1, \ldots, \lhs^e\big)$ and $\lhs^k(\vx)= \sum_{\vi{n}} \cc{\lhs}{k}{\vi{n}}\vx^\vi{n}$ for $1 \leq k \leq e$, and analogously for $\rhs$;
\item $\ve{g}=(g^1, \ldots, g^p)$ and $1 + g^k(\vx)= \sum_{\vi{n}} \cc{g}{k}{\vi{n}}\vx^\vi{n}$ for $1 \leq k \leq p$;
\item $\ve{h}=(h^1, \ldots, h^{d-(s+p)})$ and $h^k(\vx)= \sum_{\vi{n}} \cc{h}{k}{\vi{n}}\vx^\vi{n}$ for $1 \leq k \leq d-(s+p)$.
\end{itemize}

\begin{myoss}
Multi-indices $\vi{n} \in \nN^d$, although they are written as horizontal vectors, are meant to be \myemph{vertical} vectors. We shall always omit the transposition on multi-indices, but we still use subscripts to indicate their coordinates, instead of superscripts used for horizontal vectors (as in the standard notation).
\end{myoss}

We shall use the notation $\vi{n}=(\vi{n}_\vxu, \vi{n}_\vxv, \vi{n}_\vxy, \vi{n}_\vxz)$, denoting by $\ _\vxu$ the projection onto the coordinate $\vxu$ and analogously for the other coordinates, so that $\vi{n}_\vxu \in \nN^r$, $\vi{n}_\vxv \in \nN^e$, $\vi{n}_\vxy \in \nN^p$ and $\vi{n}_\vxz \in \nN^{d-(s+p)}$.

In the following, we shall need some properties of formal power series and new notations to keep the equations as compact as possible.

\begin{myoss}\label{oss:expproductrule}
Let $\vxx=(\xx^1, \ldots \xx^r) \in \nC^r$, $\mi{A} \in \nmat{a}{b}{\nN}$ and $\mi{B} \in \nmat{b}{c}{\nN}$.
By direct computation we get
$$
\big(\vxx^\mi{A}\big)^\mi{B} = \vxx^{\mi{A}\mi{B}}\mbox.
$$ 
\end{myoss}
\begin{myoss}\label{oss:powerseriesnotations}
Let $\psi:(\nC^c, 0) \rightarrow \nC^b$ be a formal map, and $\vi{i} \in \nN^b$ a multi-index. Pick $\vx=(\x^1, \ldots, \x^c)$ some coordinates at $0 \in \nC^c$. We shall need to write in formal power series expressions of the form
$$
\big(\psi(\vx)\big)^\vi{i} \in \nC[[\vx]]\mbox.
$$
Set $\vi{i}=(i_1, \ldots, i_b)$ and $\psi=(\psi^1, \ldots, \psi^b)$ with $\psi^k(\vx)=\sum_{\vi{n}} \cc{\psi}{k}{\vi{n}} \vx^{\vi{n}}$ for $k=1, \ldots, b$.
Then
\begin{equation*}
\big(\psi(\vx)\big)^\vi{i}
= \prod_{k=1}^{b} \big(\psi^k(\vx)\big)^{i_k}
= \prod_{k=1}^b \Bigg( \sum_{\vi{n}^k \in \nN^c} \cc{\psi}{k}{\vi{n}^k} \vx^{\vi{n}^k} \Bigg)^{i_k}
= \prod_{k=1}^b \prod_{l=1}^{i^k}\Bigg( \sum_{\vi{n}^{k,l} \in \nN^c} \cc{\psi}{k}{\vi{n}^{k,l}} \vx^{\vi{n}^{k,l}} \Bigg)\mbox.
\end{equation*}
Set
$$
\nindd{\vi{i}}{c}:=\big\{\vpp{N}=(\vi{n}^{1,1}, \ldots \vi{n}^{1,i^1}\ |\ \cdots\ |\ \vi{n}^{b,1}, \ldots \vi{n}^{b,i^b})\mbox{ s.t. } \vi{n}^{k,l} \in \nN^c\ \forall k,l\} \cong \nmat{c}{\abs{\vi{i}}}{\nN}\mbox.
$$
and for $\vpp{N} \in \nindd{\vi{i}}{c}$ write
$$
\psi_\vpp{N}:=\prod_{k=1}^b \prod_{l=1}^{i_k} \cc{\psi}{k}{\vi{n}^{k,l}} \in \nC \mbox, \qquad
\abs{\vpp{N}}:=\sum_{k=1}^b \sum_{l=1}^{i_k} \vi{n}^{k,l} \in \nN^c\mbox.
$$
Then we have
$$
\big(\psi(\vx)\big)^\vi{i}= \sum_{\vpp{N} \in \nindd{\vi{i}}{c}} \psi_\vpp{N} \vx^{\abs{\vpp{N}}}\mbox.
$$
When $c=d$, we shall omit the subscript and write $\nindd{\vi{i}}{d}=\nind{\vi{i}}$.
\end{myoss}

Coming back to the proof of Theorem \ref{thm:rigidprimaryresonance}, by direct computations we get
\begin{align}
\lhs^k&=\sum_{\vi{i} \in \nN^d} \left[\cc{\phi}{k}{\vi{i}}
\left(\ve{\alpha}\vxu^\mi{B}\right)^{\vi{i}_\vxu}
\big(\ve{\mu}\vxv + \ve{\rho}(\vx)\big)^{\vi{i}_\vxv}
\Big(\ve{\beta}\vxu^\mi{C}\vxy^\mi{D} \big(\one +  \ve{g}(\vx)\big)\Big)^{\vi{i}_\vxy}
\big(\ve{h}(\vx)\big)^{\vi{i}_\vxz}
\right]\nonumber
\\
&=\sum_{\vi{i} \in \nN^d} \cc{\phi}{k}{\vi{i}}
\ve{\alpha}^{\vi{i}_\vxu} \vxu^{\mi{B} \vi{i}_\vxu}
\ve{\beta}^{\vi{i}_\vxy} \vxu^{\mi{C} \vi{i}_\vxy} \vxy^{\mi{D}\vi{i}_\vxy}
\left(
\sum_{\vpp{I} \in \nind{\vi{i}_\vxv}} \rho_\vpp{I} \vx^{\abs{\vpp{I}}}
\sum_{\vpp{J} \in \nind{\vi{i}_\vxy}} g_\vpp{J} \vx^{\abs{\vpp{J}}}
\sum_{\vpp{K} \in \nind{\vi{i}_\vxz}} h_\vpp{K} \vx^{\abs{\vpp{K}}} 
\right)\label{eqn:rprlhsexpr}
\\ \nonumber \\
\rhs^k&=\sum_{\vi{j} \in \nN^{r+e}} \cc{\widetilde{\rho}}{k}{\vi{j}}
\vxu^{\vi{j}_\vxu}
\big(\ve{\phi}(\vx)\big)^{\vi{j}_\vxv}
=\sum_{\vi{j} \in \nN^{r+e}} \cc{\widetilde{\rho}}{k}{\vi{j}}
\vxu^{\vi{j}_\vxu}
\left(\sum_{\vpp{H} \in \nind{\vi{j}_\vxv}} \phi_\vpp{H} \vx^{\abs{\vpp{H}}}\right)\mbox, \label{eqn:rprrhsexpr}
\end{align}
for $k=1, \ldots, e$.

Expressing explicitly the coefficients of $\lhs^k$ and $\rhs^k$ written in formal power series, from \eqref{eqn:rprlhsexpr} and \eqref{eqn:rprrhsexpr} respectively we obtain:
$$
\cc{\lhs}{k}{\vi{n}} = 
\sumaa{\substack{\vi{i}\in \nN^d \\ \vpp{I} \in \nind{\vi{i}_\vxv}, \vpp{J} \in \nind{\vi{i}_\vxy}, \vpp{K} \in \nind{\vi{i}_\vxz}\\
\cond{1}}}{1.2cm}
\cc{\phi}{k}{\vi{i}} \ve{\alpha}^{\vi{i}_\vxu}\ve{\beta}^{\vi{i}_\vxy} \rho_{\vpp{I}} g_\vpp{J} h_{\vpp{K}}
\mbox, 
\qquad
\cc{\rhs}{k}{\vi{n}} =
\sumaa{\substack{\vi{j} \in \nN^{r+e}
\\ \vpp{H} \in \nind{\vi{j}_\vxv} \\ \cond{2}}}{0.3cm}
\cc{\widetilde{\rho}}{k}{\vi{j}}\phi_{\vpp{H}} 
\mbox,
$$
for $k=1, \ldots, e$ and $\vi{n} \in \nN^d$; moreover
$$
\cond{1}=
\left\{
\begin{array}{l}
\mi{B} \vi{i}_\vxu + \mi{C} \vi{i}_\vxy + \abs{\vpp{I}}_\vxu +\abs{\vpp{J}}_\vxu +\abs{\vpp{K}}_\vxu = \vi{n}_\vxu \\
\abs{\vpp{I}}_{\vxv} +\abs{\vpp{J}}_{\vxv} +\abs{\vpp{K}}_{\vxv}= \vi{n}_{\vxv} \\
\mi{D} \vi{i}_\vxy + \abs{\vpp{I}}_\vxy +\abs{\vpp{J}}_\vxy +\abs{\vpp{K}}_\vxy= \vi{n}_\vxy \\
\abs{\vpp{I}}_{\vxz} +\abs{\vpp{J}}_{\vxz} +\abs{\vpp{K}}_{\vxz}= \vi{n}_{\vxz}
\end{array}
\right.
\mbox,
$$
and
$$
\cond{2}=
\left\{
\begin{array}{l}
\vi{j}_\vxu + \abs{\vpp{H}}_\vxu = \vi{n}_\vxu \\
\abs{\vpp{H}}_{\vxv} = \vi{n}_{\vxv} \\
\abs{\vpp{H}}_{\vxy} = \vi{n}_{\vxy} \\
\abs{\vpp{H}}_{\vxz} = \vi{n}_{\vxz}
\end{array}
\right.
\mbox.
$$

We want to solve the equation 
\begin{equation}\label{eqn:primarycoefftosolve}
\cc{\eq}{k}{\vi{n}}:= \cc{\rhs}{k}{\vi{n}}-\cc{\lhs}{k}{\vi{n}}=0
\end{equation}
for every $k$ and $\vi{n}$, where the unknowns are the coefficients $\cc{\phi}{k}{\vi{n}}$ of $\ve{\phi}$ and $\cc{\widetilde{\rho}}{k}{\vi{n}}$ of $\ve{\widetilde{\rho}}$.

To understand the combinatorics of \eqref{eqn:primarycoefftosolve}, we need a partial order and a total order on indices in $\nN^d$.
Set $\vi{n}=(n_1, \ldots, n_d)$ and $\vi{m}=(m_1, \ldots, m_d)$.
\begin{trivlist}
\item[Partial order $\preceq$:] we say that $\vi{m} \preceq \vi{n}$ iff we have $m_k \leq n_k$ for every $k=1, \ldots, d$.
\item[Total order $\leq$:] we say that $\vi{m} \leq \vi{n}$ iff $(\abs{\vi{m}}, m_1, \ldots, m_d) \leq_{\on{lex}} (\abs{\vi{n}}, n_1, \ldots, n_d)$, where $\leq_{\on{lex}}$ is the lexicographic order (on $\nN^{d+1}$).

For example, for $d=3$ we have:
\begin{align*}
&(0,0,0)<\\
&(0,0,1)<(0,1,0)<(1,0,0)<\\
&(0,0,2)<(0,1,1)<(0,2,0)<(1,0,1)<(1,1,0)<(2,0,0)<\\
&(0,0,3)<(0,1,2)<\cdots<(2,1,0)<(3,0,0)<\\
&\hspace{0.6cm}\vdots
\end{align*}
\end{trivlist}
We notice that if $\vi{m} \prec \vi{n}$ then $\vi{m} < \vi{n}$. Moreover if $\vi{m}^\prime \leq \vi{n}^\prime$ and $\vi{m}^\second \leq \vi{n}^\second$ then $\vi{m}^\prime + \vi{m}^\second \leq \vi{n}^\prime + \vi{n}^\second$.

\begin{mylem}\label{lem:computationpsiH}
Let $\ve{\psi}:(\nC^d, 0) \rightarrow (\nC^b, 0)$, $\vi{j} \in \nN^b$ and $\vpp{H} \in \nind{\vi{j}}$.
Take coordinates $\vx=(\x^1, \ldots, \x^d) \in \nC^d$, and set $\psi=(\psi^1, \ldots, \psi^b)$ with $\psi^k(\vx)=\sum_{\vi{n}}\cc{\psi}{k}{\vi{n}} \vx^{\vi{n}}$ for every $k=1, \ldots, b$.

For any $k=1, \ldots, d$, let $\vi{e}^k \in \nN^d$ be the multi-index with $1$ in the $k$-th coordinate and $0$ in all the others.
Suppose there exists $0 \leq c \leq d-b-1$ such that $\cc{\psi}{k}{\vi{n}} = 0$ for every $\vi{n} < \vi{e}_{c+k}$, $k=1, \ldots, b$.

Then $\psi_\vpp{H} = 0$ for $\abs{\vpp{H}} < (\vi{0}_{c}, \vi{j}, \vi{0}_{d-c-b})$ (where $\vi{0}_c \in \nN^c$ and $\vi{0}_{d-c-b} \in \nN^{d-c-b}$).
Moreover:
\begin{enumerate}[(i)]
\item if $\cc{\psi}{k}{\vi{e}^{c+k}}=0$ for $k=1, \ldots, b$, then $\psi_\vpp{H}=0$ for $\abs{\vpp{H}}\leq(\vi{0}_{c}, \vi{j}, \vi{0}_{d-c-b})$ if $\vi{j} \neq 0$;
\item if $\psi^k(\vx)=\zeta^k \x^{c+k} +\ \on{h.o.t.}$ for $k=1, \ldots, b$, then $\psi_\vpp{H} \neq 0$ only if one of the following conditions is satisfied:
\begin{itemize}
\item $\abs{\vpp{H}}=(\vi{0}_{c}, \vi{j}, \vi{0}_{d-c-b})$, and in this case $\vpp{H}$ is uniquely determined in $\nind{\vi{j}}$ and $\psi_\vpp{H}=\zeta^{\vi{j}}$, where $\zeta=(\zeta^1, \ldots, \zeta^b)$;
\item $\abs{\abs{\vpp{H}}} > \abs{\vi{j}}$, where $\abs{\abs{\vpp{H}}} \in \nN$ denotes the sum of all elements of $\abs{\vpp{H}}\in \nN^d$, and hence the sum of all elements of $\vpp{H} \in \nmat{d}{\abs{\vi{j}}}{\nN}$.
\end{itemize}
\end{enumerate}
\end{mylem}
\begin{mydim}
Set $\vi{j}=(j_1, \ldots, j_b)$, and write explicitly
$$
\vpp{H}=
\begin{mybmatrix}{3}
\vi{h}^{1,1}\cdots \vi{h}^{1,j_1} & \vi{h}^{2,1} \cdots \vi{h}^{2,j_2} & \cdots & \vi{h}^{b,1} \cdots \vi{h}^{b,j_b}
\end{mybmatrix}\mbox,
$$
where $\vi{h}^{k,l} \in \nN^d$ is a multi-index for every $k=1, \ldots, b$ and $l=1, \ldots, j_k$.

To have $\psi_\vpp{H} \neq 0$, we must have $\cc{\psi}{k}{\vi{h}^{k,l}} \neq 0$ for every $k$ and $l$.

Thanks to our assumption, $\cc{\psi}{k}{\vi{h}^{k,l}} \neq 0$ only if $\vi{h}^{k,l} \geq \vi{e}^{c+k}$.
Then we have that $\phi_\vpp{H} \neq 0$ only if
\begin{equation}\label{eqn:indexestimate}
\abs{\vpp{H}} = \sum_{k=1}^b\sum_{l=1}^{j_k} \vi{h}^{k,l} \geq \sum_{k=1}^b\sum_{l=1}^{j_k} \vi{e}^{c+k} = \sum_{k=1}^b j_k \vi{e}^{c+k} = (\vi{0}_c, \vi{j}, \vi{0}_{d-c-b})\mbox.  
\end{equation}

\begin{enumerate}[(i)]
\item Assume that $\cc{\psi}{k}{\vi{h}^{k,l}} \neq 0$ only if $\vi{h}^{k,l} > \vi{e}^{c+k}$.
Since $\vi{j} \neq \vi{0}$, the sums in \eqref{eqn:indexestimate} are not empty, and the inequality is strict.
\item Since $\phi^k(\vx)- \zeta^k \x^{c+k}$ is at least of order $2$, $\cc{\psi}{k}{\vi{n}}\neq 0$ only if $\vi{n}=\vi{e}^{c+k}$ or $\abs{\vi{n}} \geq 2$.

\begin{itemize}
\item If $\vi{h}^{k,l} = \vi{e}^{c+k}$ for every $k, l$, then
$$
\abs{\vpp{H}}=
\abs{
\begin{mybmatrix}{3}
0 \cdots 0 & 0 \cdots 0 & \cdots & 0 \cdots 0 \\
\vdots & \vdots & \vdots & \vdots \\
0 \cdots 0 & 0 \cdots 0 & \cdots & 0 \cdots 0 \\
\hline
1 \cdots 1 & 0 \cdots 0 & \cdots & 0 \cdots 0 \\
0 \cdots 0 & 1 \cdots 1 & \cdots & 0 \cdots 0 \\
\vdots & \vdots & \vdots & \vdots \\
0 \cdots 0 & 0 \cdots 0 & \cdots & 1 \cdots 1 \\
\hline
0 \cdots 0 & 0 \cdots 0 & \cdots & 0 \cdots 0 \\
\vdots & \vdots & \vdots & \vdots \\
0 \cdots 0 & 0 \cdots 0 & \cdots & 0 \cdots 0
\end{mybmatrix}
}
=
\begin{mymatrix}{1}
0\\
\vdots\\
0\\
\hline
j_1\\
j_2\\
\vdots\\
j_b\\
\hline
0\\
\vdots\\
0
\end{mymatrix}
=(\vi{0}_c,\vi{j}, \vi{0}_{d-c-b})
\mbox.
$$
In this case,
$$
\phi_\vpp{H}=\prod_{k=1}^b \prod_{l=1}^{\vi{j}_k} \cc{\psi}{k}{\vi{e}^{c+k}}= \prod_{k=1}^b (\zeta^k)^{j_k} = \zeta^{\vi{j}} \mbox.
$$
\item If there exist $k$ and $l$ such that $\abs{\vi{h}^{k,l}} \geq 2$, then to have $\phi_\vpp{H} \neq 0$ we must have 
$$
\abs{\abs{\vpp{H}}} = \sum_{k=1}^e \sum_{l=1}^{\vi{j}_k} \abs{\vi{h}^{k,l}} > \sum_{k=1}^e \vi{j}_k = \abs{\vi{j}}\mbox.
$$
\end{itemize}
\end{enumerate}
\end{mydim}

Recall that $e$ is the number of components of $\vxv$, and hence of $\lhs$ and $\rhs$.
We shall need a weight on indices $(k, \vi{n}) \in \{1, \ldots, e\} \times \nN^d$.

\begin{mydef}\label{def:primaryweight}
Let $k \in \{1, \ldots, e\}$ be an integer and $\vi{n}$ be a multi-index (in $\nN^d$ or $\nN^s$).
We call \mydefin{weight} of $(k, \vi{n})$ the value
$$
\wa{k}{\vi{n}}=\wa{k}{\abs{\vi{n}}} :=\abs{\vi{n}} + \frac{k}{e} \in \frac{\nN}{e}\mbox.
$$
\end{mydef}

Notice that for every $W \in \nN/e$, there are only finitely many $(k, \vi{n})$ such that $\wa{k}{\vi{n}} \leq W$.

\begin{mylem}\label{lem:primaryrhslot}
We have
\begin{equation*}
\cc{\rhs}{k}{\vi{n}} = \kronecker{\vi{n}_\vxy}{\vi{0}} \kronecker{\vi{n}_\vxz}{\vi{0}} \cc{\widetilde{\rho}}{k}{(\vi{n}_\vxu, \vi{n}_\vxv)} + \mu_k \cc{\phi}{k}{\vi{n}} 
+ Q_{k, \abs{\vi{n}}}(\cc{\phi}{l}{\vi{m}}, \cc{\widetilde{\rho}}{l}{(\vi{m}_\vxu, \vi{m}_\vxv)})\mbox,
\end{equation*}
where $\kronecker{}{}$ denotes the Kronecker's delta function, and $Q_{k, \abs{\vi{n}}}$ is a polynomial in the variables $\cc{\phi}{l}{\vi{m}}$ and $\cc{\widetilde{\rho}}{l}{(\vi{m}_\vxu, \vi{m}_\vxv)}$ satisfying
$$
\wa{l}{\vi{m}} < \wa{k}{\abs{\vi{n}}}\mbox.
$$
\end{mylem}
In order to simplify notations, we shall simply write
\begin{equation}\label{eqn:primaryrhslot}
\cc{\rhs}{k}{\vi{n}} = \kronecker{\vi{n}_\vxy}{\vi{0}} \kronecker{\vi{n}_\vxz}{\vi{0}} \cc{\widetilde{\rho}}{k}{(\vi{n}_\vxu, \vi{n}_\vxv)} + \mu_k \cc{\phi}{k}{\vi{n}} 
+ \lllot{k}{\abs{\vi{n}}}{\phi, \widetilde{\rho}}\mbox,
\end{equation}
where $\lllot{k}{\abs{\vi{n}}}{\phi, \widetilde{\rho}}$ stands for a suitable polynomial in the variables $\cc{\phi}{l}{\vi{m}}$ and $\cc{\widetilde{\rho}}{l}{(\vi{m}_\vxu, \vi{m}_\vxv)}$ satisfying $\wa{l}{\vi{m}} < \wa{k}{\abs{\vi{n}}}$.
We shall also omit $\widetilde{\rho}$ when the polynomial does not depend on any coefficient $\cc{\widetilde{\rho}}{l}{(\vi{m}_\vxu, \vi{m}_\vxv)}$.
\begin{mydim}
Set $W:=\wa{k}{\abs{n}}$.
From the first equation of $\cond{2}$, we get $\vi{j}_\vxu \preceq \vi{n}_\vxu$, and in particular $\abs{\vi{j}_\vxu} \leq \abs{\vi{n}_\vxu}$.
From Lemma \ref{lem:computationpsiH}.(ii) we can have two cases when $\phi_\vpp{H} \neq 0$.
\begin{itemize}
\item Either $\vi{j}_\vxv = \abs{\vpp{H}}_\vxv = \vi{n}_\vxv$, and in this case the term $\cc{\widetilde{\rho}}{k}{\vi{j}}$ with the biggest weight is given by $\vi{j}_\vxu=\vi{n}_\vxu$. Its weight is $\leq W$, and the equality holds only if $\vi{n}_\vxy=\vi{0}$ and $\vi{n}_\vxz=\vi{0}$, when we get the first term of \eqref{eqn:primaryrhslot}.
\item Or
$$
\abs{\vi{j}} = \abs{\vi{j}_\vxu} + \abs{\vi{j}_\vxv} < \abs{\vi{j}_\vxu} + \abs{\abs{\vpp{H}}} = \abs{\vi{j}_\vxu} + \abs{\vi{n}}-\abs{\vi{j}_\vxu} = \abs{\vi{n}}\mbox,
$$
and in this case the weight strictly less than $W$.
\end{itemize}

Still from $\cond{2}$, we get $\abs{H} \prec \vi{n}$.
It follows that the only way to have $\phi_\vpp{H} \neq 0$ and with some $\cc{\phi}{l}{\vi{m}}$ with $\wa{l}{m}\geq W$ is to have $\vpp{H}$ made by just a column in position $l \geq k$, given by $\vi{n}$.
In this case we get $\vi{j}_\vxv=\vi{e}^l$, $\phi_\vpp{H}=\cc{\phi}{l}{\vi{n}}$, and from the first equation of $\cond{2}$ we get $\vi{j}_\vxu=\vi{0}$.
Since $\cc{\widetilde{\rho}}{k}{\vi{e}^l}=0$ for $l>k$ and $\cc{\widetilde{\rho}}{k}{\vi{e}^k}=\mu^k$, we get the second term of \eqref{eqn:primaryrhslot}.
\end{mydim}

\begin{mylem}\label{lem:primarylhslot}
We have
\begin{equation}\label{eqn:primarylhslot}
\cc{\lhs}{k}{\vi{n}} = \kronecker{\vi{n}_\vxy}{\vi{0}} \kronecker{\vi{n}_\vxz}{\vi{0}} \ve{\alpha}^{\mi{B}^{-1} \vi{n}_\vxu} \ve{\mu}^{\vi {n}_\vxv} \cc{\phi}{k}{(\mi{B}^{-1}\vi{n}_\vxu, \vi{n}_\vxv, \vi{0}, \vi{0})} + \lllot{k}{\abs{\vi{n}}}{\phi}\mbox.
\end{equation}
\end{mylem}
\begin{mydim}
Set $W:=\wa{k}{\abs{n}}$.
Thanks to Lemma \ref{lem:computationpsiH} we get that $\rho_\vpp{I} \neq 0$ only if $\abs{\vpp{I}} \geq (\vi{0},\vi{i}_\vxv,\vi{0},\vi{0})$.
Thanks to Lemma \ref{lem:computationpsiH}.(i) we get that $h_\vpp{K} \neq 0$ only if $\abs{\vpp{K}} > (\vi{0},\vi{0},\vi{0},\vi{i}_\vxz)$ when $\vi{i}_\vxz \neq \vi{0}$.

From the first equation in $\cond{1}$ we get $\mi{B} \vi{i}_\vxu \preceq \vi{n}_\vxu$, hence $\vi{i}_\vxu \preceq \mi{B}^{-1} \vi{n}_\vxu$ and in particular $\abs{\vi{i}_\vxu} \leq \abs{\vi{n}_\vxu}$.
Notice that the equality on modules holds only if $\vi{i}_\vxu=\mi{B}^{-1} \vi{n}_\vxu$.
From the third equation we get $\abs{\vi{i}_\vxy} < \abs{\mi{D} \vi{i}_\vxy} \leq \abs{\vi{n}_\vxy}$ for $\vi{i}_\vxy \neq 0$.

Then we get
\begin{align*}
\abs{\vi{i}}=&\abs{\vi{i}_\vxu}+\abs{\vi{i}_\vxv}+\abs{\vi{i}_\vxy}+\abs{\vi{i}_\vxz}\\
\leq& \abs{\vi{i}_\vxu}+\abs{\abs{\vpp{I}}}+\abs{\vi{i}_\vxy}+\abs{\abs{\vpp{K}}} \\
=&\abs{\vi{n}} - \abs{\abs{\vpp{J}}} - \abs{\mi{C} \vi{i}_\vxy} - \big(\abs{\mi{D} \vi{i}_\vxy} - \abs{\vi{i}_\vxy}\big)\\
\leq \abs{\vi{n}}\mbox,
\end{align*}
where the equality can hold only if $\vi{i}_\vxy=0$ and $\vi{i}_\vxz=0$.

It follows that terms $\cc{\phi}{k}{\vi{i}}$ such that $\wa{k}{\vi{i}} \geq W$ appear only when $\vi{i}_\vxy=\vi{n}_\vxy=\vi{0}$ and $\vi{i}_\vxz=\vi{n}_\vxz=\vi{0}$.
In this case, $\vpp{J}=\vpp{K}=\emptyset$, and $g_\emptyset=h_\emptyset=1$.
The third equation of $\cond{1}$ gives $\abs{I}_\vxy=\vi{0}$, while the fourth gives $\abs{I}_\vxz=\vi{0}$.
The second equation of $\cond{1}$ gives $\abs{I}_\vxv=\vi{n}_\vxv$.
To have a term $\cc{\phi}{k}{\vi{i}}$ of weight $W$, we need to have then $\vi{i}_\vxu=\mi{B}^{-1} \vi{n}_\vxu$, from which it follows $\abs{I}_\vxu=0$.
Then the second equation of $\cond{1}$ gives $\vi{i}_\vxv \leq \abs{I}_\vxv=\vi{n}_\vxv$.
But Lemma \ref{lem:computationpsiH} says that $\rho_\vpp{I} = 0$ for $\abs{I}_\vxv > \vi{i}_\vxv$.
Hence the only term that appears is for $\vi{i}_\vxv=\vi{n}_\vxv$.
Following the computation of the proof of Lemma \ref{lem:computationpsiH}.(ii), we get $\rho_\vpp{I}=\mu^{\vi{n}_\vxv}$, and the statement.
\end{mydim}

Thanks to Lemmas \ref{lem:primaryrhslot} and \ref{lem:primarylhslot}, $\cc{\eq}{k}{\vi{n}}=0$ becomes
\begin{equation}\label{eqn:primarylot}
\mu^k \cc{\phi}{k}{\vi{n}} + \kronecker{\vi{n}_\vxy}{\vi{0}} \kronecker{\vi{n}_\vxz}{\vi{0}} \left( \cc{\widetilde{\rho}}{k}{(\vi{n}_\vxu, \vi{n}_\vxv)} -  \ve{\alpha}^{\mi{B}^{-1} \vi{n}_\vxu} \ve{\mu}^{\vi {n}_\vxv} \cc{\phi}{k}{(\mi{B}^{-1} \vi{n}_\vxu, \vi{n}_\vxv, \vi{0}, \vi{0})} \right) = \lllot{k}{\abs{n}}{\phi, \widetilde{\rho}}\mbox.
\end{equation}
This affine equation, where the unknowns are $\cc{\phi}{k}{\vi{n}}$ and $\cc{\widetilde{\rho}}{k}{(\vi{n}_\vxu, \vi{n}_\vxv)}$, has always a solution.
At this point we conjugated $f$ to a map $\widetilde{f}$ as in \eqref{eqn:rigidsplitprimaryresonance}, but with $\widetilde{\rho}:(\nC^s,0) \rightarrow (\nC^e,0)$ a (vector of) formal power series.
Next we show that we can solve the conjugacy relation \eqref{eqn:primarycoefftosolve} and get $\widetilde{\rho}$ polynomial with only primary resonant monomials.

We solve $\cc{\eq}{k}{\vi{n}}=0$ inductively on $\wa{k}{\vi{n}}$ as follows.
\begin{trivlist}
\item For $\wa{k}{\vi{n}} \leq 2$, i.e., if $\abs{\vi{n}} \leq 1$, we set $\cc{\phi}{k}{\vi{n}}:=1$ if $\vi{n}=\vi{e}^{r+k}$ and $0$ otherwise, while $\cc{\widetilde{\rho}}{k}{(\vi{n}_\vxu, \vi{n}_\vxv)}:= \cc{\rho}{k}{(\vi{n}_\vxu, \vi{n}_\vxv, \vi{0}, \vi{0})}$. An easy computation shows that $\cc{\eq}{k}{\vi{n}}=0$ holds for these values.
\item Set $2<W \in \nN/e$, and suppose that we have determined $\cc{\phi}{l}{\vi{m}}$ and $\cc{\widetilde{\rho}}{l}{(\vi{m}_\vxu, \vi{m}_\vxv)}$ for $\wa{l}{\vi{m}} < W$ satisfying $\cc{\eq}{l}{\vi{m}}=0$ when $\wa{l}{\vi{m}}<W$.
We want to solve \eqref{eqn:primarylot} for $\wa{k}{\vi{n}}=W$.

Notice that $\lllot{k}{\abs{\vi{n}}}{\phi, \widetilde{\rho}}$ is a polynomial that depend on $\cc{\phi}{l}{\vi{m}}$ and $\cc{\widetilde{\rho}}{l}{(\vi{m}_\vxu, \vi{m}_\vxv)}$ only for weights strictly less than $W$.
Hence thanks to the induction hypothesis, $\lllot{k}{\abs{\vi{n}}}{\phi, \widetilde{\rho}}$ is a known value in $\nC$.
\begin{enumerate}[$1$)]
\item Suppose $(\vi{n}_\vxy, \vi{n}_\vxz) \neq (\vi{0}, \vi{0})$. Then \eqref{eqn:primarylot} becomes
$$
\mu^k \cc{\phi}{k}{\vi{n}} = \lllot{k}{\abs{n}}{\phi, \widetilde{\rho}}\mbox,
$$
and there exists a unique $\cc{\phi}{k}{\vi{n}}$ that solves the equation.
\item Suppose $(\vi{n}_\vxy, \vi{n}_\vxz) = (\vi{0}, \vi{0})$. Then \eqref{eqn:primarylot} becomes
\begin{equation}\label{eqn:primarylot00}
-\mu^k \cc{\phi}{k}{\vi{n}} + \ve{\alpha}^{\mi{B}^{-1} \vi{n}_\vxu} \ve{\mu}^{\vi {n}_\vxv} \cc{\phi}{k}{(\mi{B}^{-1} \vi{n}_\vxu, \vi{n}_\vxv, \vi{0}, \vi{0})} = \cc{\widetilde{\rho}}{k}{(\vi{n}_\vxu, \vi{n}_\vxv)}  + \lllot{k}{\abs{n}}{\phi, \widetilde{\rho}}\mbox.
\end{equation}
\begin{enumerate}[$2.1$)]
\item Suppose that $\vi{n}_\vxu=\mi{B}^{-1} \vi{n}_\vxu$: we have two cases.

Suppose $\mu^k \neq \ve{\alpha}^{\vi{n}_\vxu} \ve{\mu}^{\vi {n}_\vxv}$, i.e., $\vxu^{\vi{n}_\vxu} \vxv^{\vi{n}_\vxv}$ is not primary resonant for the $k$-th coordinate.
Then we can put $\cc{\widetilde{\rho}}{k}{(\vi{n}_\vxu, \vi{n}_\vxv)}=0$ and there exists a unique $\cc{\phi}{k}{\vi{n}}$ that solves the equation.

Suppose $\mu^k = \ve{\alpha}^{\vi{n}_\vxu} \ve{\mu}^{\vi {n}_\vxv}$, i.e., $\vxu^{\vi{n}_\vxu} \vxv^{\vi{n}_\vxv}$ is primary resonant for the $k$-th coordinate.
Then \eqref{eqn:primarylot00} does not depend on $\cc{\phi}{k}{\vi{n}}$ (we put it equal to $0$), and there exists a unique $\cc{\widetilde{\rho}}{k}{(\vi{n}_\vxu, \vi{n}_\vxv)}$ that solves the equation.

\item In the general case, let $\widetilde{\eta}$ be the smallest number in $\nN^*$ such that $\vi{n}_\vxu=\mi{B}^{\widetilde{\eta}} \vi{n}_\vxu$.
Set $\vi{n}_\vxu^{(l)}:=\mi{B}^{-l} \vi{n}_\vxu$ for $l=0, \ldots, \widetilde{\eta}-1$.
We consider the equation \eqref{eqn:primarylot00} for $\vi{n}_\vxu, \vi{n}_\vxu^{(1)}, \ldots, \vi{n}_\vxu^{(1-\widetilde{\eta})}$ simultaneously (while we fix $\vi{n}_\vxv$).
In this case we get the following linear system:
\begin{equation}\label{eqn:primarylinearsystem}
\begin{mymatrix}{1}
\!\!\!\cc{\phi}{k}{(\vi{n}_\vxu, \vi{n}_\vxv, \vi{0}, \vi{0})}\!\!\!\\
\!\!\!\cc{\phi}{k}{(\vi{n}_\vxu^{(1)}, \vi{n}_\vxv, \vi{0}, \vi{0})}\!\!\!\\
\vdots \\
\vdots \\
\!\!\!\cc{\phi}{k}{(\vi{n}_\vxu^{(\widetilde{\eta}-1)}, \vi{n}_\vxv, \vi{0}, \vi{0})}\!\!\!
\end{mymatrix}^T
\begin{mymatrix}{5}
-\mu^k & 0 & \cdots & 0 & \Big.\!\!\!\ve{\alpha}^{\vi{n}_\vxu} \ve{\mu}^{\vi{n}_\vxv}\!\!\!\Big. \\
\!\!\!\ve{\alpha}^{\vi{n}_\vxu^{(1)}} \ve{\mu}^{\vi{n}_\vxv} \!\!\! & -\mu^k & \ddots & \ddots & 0 \\
0 & \!\!\! \ve{\alpha}^{\vi{n}_\vxu^{(2)}} \ve{\mu}^{\vi{n}_\vxv} \!\!\! & \ddots & \ddots & \vdots \\
\vdots & \ddots & \ddots & -\mu^k & 0 \\
0 & \cdots & 0 & \!\!\! \ve{\alpha}^{\vi{n}_\vxu^{(\widetilde{\eta}-1)}} \ve{\mu}^{\vi{n}_\vxv} \!\!\! & -\mu^k
\end{mymatrix}
=
\begin{mymatrix}{1}
\!\!\!\cc{\widetilde{\rho}}{k}{(\vi{n}_\vxu, \vi{n}_\vxv)}\!\!\!\\
\!\!\!\cc{\widetilde{\rho}}{k}{(\vi{n}_\vxu^{(1)}, \vi{n}_\vxv)}\!\!\!\\
\vdots \\
\vdots \\
\!\!\!\cc{\widetilde{\rho}}{k}{(\vi{n}_\vxu^{(\widetilde{\eta}-1)}, \vi{n}_\vxv)}\!\!\!
\end{mymatrix}^T
+ \lot\mbox,
\end{equation}
where
$$
\lot=\lllot{k}{\abs{\vi{n}}}{\phi, \widetilde{\rho}}\mbox.
$$
for each coordinate, since $\wa{k}{(\vi{n}_\vxu^{(l)},\vi{n}_\vxv,\vi{0},\vi{0})}=W$ for each $l=0, \ldots, \widetilde{\eta}-1$.

By direct computation, the determinant of the matrix in \eqref{eqn:primarylinearsystem} is given (up to sign) by
$$
\left(\prod_{l=0}^{\widetilde{\eta}-1} \ve{\alpha}^{\vi{n}_\vxu^{(l)}} \right) \ve{\mu}^{\widetilde{\eta} \vi{n}_\vxv} - (\mu^k)^{\widetilde{\eta}}\mbox.
$$
Recall that the $\eta$ in the definition of primary resonances \eqref{eqn:defprimaryresonance} is the order of $\mi{B}$.
In particular, we have $\widetilde{\eta} \mid \eta$.
It follows that the linear system \eqref{eqn:primarylinearsystem} is invertible iff $(\vi{n}_\vxu, \vi{n}_\vxv)$ is not primary resonant.
Notice also that $(\vi{n}_\vxu, \vi{n}_\vxv)$ is primary resonant iff $(\vi{n}_\vxu^{(l)},\vi{n}_\vxv)$ is for every $l=0, \ldots, \widetilde{\eta}-1$.

If $(\vi{n}_\vxu, \vi{n}_\vxv)$ is not primary resonant, we can put $\cc{\widetilde{\rho}}{k}{(\vi{n}_\vxu^{(l)}, \vi{n}_\vxv)}=0$ for every $l=0, \ldots, \widetilde{\eta}-1$ and there exists a unique $(\cc{\phi}{k}{(\vi{n}_\vxu^{(l)}, \vi{n}_\vxv, \vi{0}, \vi{0})})$ for $l=0, \ldots, \widetilde{\eta}-1$ that solves the linear system \eqref{eqn:primarylinearsystem}.

If $(\vi{n}_\vxu, \vi{n}_\vxv)$ is primary resonant, we can put $\cc{\phi}{k}{(\vi{n}_\vxu^{(l)}, \vi{n}_\vxv, \vi{0}, \vi{0})}=0$ for every $l=0, \ldots, \widetilde{\eta}-1$ and there exists a unique $(\cc{\widetilde{\rho}}{k}{(\vi{n}_\vxu^{(l)}, \vi{n}_\vxv)})$ for $l=0, \ldots, \widetilde{\eta}-1$ that solves the linear system \eqref{eqn:primarylinearsystem}. 
\end{enumerate}
\end{enumerate}
\end{trivlist}
We have defined the conjugation $\Phi$ as an invertible formal map: we can then define $\widetilde{g}$ and $\widetilde{h}$ such that the conjugacy relation holds for all coordinates.
\end{mystep}

\begin{mystep}[Step $2$]

The proof of the convergence of the conjugacy map is completely analogous to the proof of the Poincar\'e-Dulac theorem (see, e.g., \cite{sternberg:localcontractionstheorempoincare}, \cite{rosay-rudin:holomorphicmaps} or \cite[Chapter $4$]{berteloot:methodeschangementechelles}).

Pick $0 < \Lambda < 1$ such that $\Lambda > \on{specrad}(df_0)$ the spectral radius of the differential $df_0$ of $f$ at $0$, and take $N$ such that $\Lambda^N<\abs{\mu_k}$ for every $k=1, \ldots, e$.

For proving the formal result, we introduced a weight, and noticed that for every $W \in \nN/e$, there are only finitely many $(k, \vi{n})$ such that $\wa{k}{\vi{n}} \leq W$.
It follows that there exist $M >0$ and a polynomial (hence holomorphic) change of coordinates that conjugates $\ve{f}$ with a map of the form
\begin{equation}\label{eqn:rigidsplitprimaryresonancewithrest}
(\vxu, \vxv, \vxy, \vxz) \mapsto \Big(\ve{\alpha} \vxu^{\mi{B}}, \ve{\mu} \vxv + \ve{\rho}(\vxu,\vxv) + \ve{R}(\vxu, \vxv, \vxy, \vxz), \ve{\beta} \vxu^\mi{C}\vxy^\mi{D} \big(\one +  \ve{g}(\vxu,\vxv,\vxy,\vxz)\big), \ve{h}(\vxu,\vxv,\vxy,\vxz) \Big)\mbox,
\end{equation}
with the same conditions as for \eqref{eqn:rigidsplitprimaryresonance}, and $\ve{R}:(\nC^d,0) \rightarrow (\nC^e, 0)$ such that
$$
\norm{\ve{R}(\vx)} \leq M \norm{\vx}^N
$$
for $\norm{\vx}$ small enough.

Notice that there are no primary resonances $\vxu^{\vi{n}_\vxu}\vxv^{\vi{n}_\vxv}$ such that $\abs{\vi{n}_\vxu} + \abs{\vi{n}_\vxv} \geq N$. 
Indeed, suppose $\vxu^{\vi{n}_\vxu}\vxv^{\vi{n}_\vxv}$ is resonant for the $k$-th coordinate for a suitable $1 \leq k \leq e$. Then from \eqref{eqn:defprimaryresonance} we would have
$$
\abs{(\mu^k)^\eta} = \abs{\ve{\xi}^{\vi{n}_\vxu}\ve{\mu}^{\eta \vi{n}_\vxv}} < \Lambda^{\eta N} < \abs{\mu^k}^\eta\mbox,
$$
that gives a contradiction.

Set $\ve{R}=(R^1, \ldots, R^e)$.
We now proceed by induction on $k=1, \ldots, e$ and prove that we can conjugate $f$ with a germ of the form \eqref{eqn:rigidsplitprimaryresonancewithrest}, with $R^l \equiv 0$ for any $l\leq k$. 
If $k=0$, there is nothing to prove.
Suppose that $\ve{f}$ is of the form \eqref{eqn:rigidsplitprimaryresonancewithrest}, with $R^l\equiv 0$ for $l < k$.
The induction step will consist in proving that we can conjugate $\ve{f}$ with $\ve{\widetilde{f}}:(\nC^d, 0) \rightarrow (\nC^d, 0)$ of the form \eqref{eqn:rigidsplitprimaryresonancewithrest}, with $\ve{\widetilde{R}}=(\widetilde{R}^1, \ldots, \widetilde{R}^e)$ instead of $\ve{R}$ such that $\widetilde{R}^l \equiv 0$ for $l \leq k$.

Consider a local diffeomorphism $\Phi:(\nC^d, 0) \rightarrow (\nC^d, 0)$ of the form
$$
\Phi(\vx)=(\vxu, \xv^1, \ldots, \xv^{k-1}, \xv^k + \phi^k(\vx), \xv^{k+1}, \ldots, \xv^e, \vxy, \vxz)\mbox,
$$
where $\phi^k:(\nC^d, 0) \rightarrow (\nC, 0)$ is of order at least $2$.

Thanks to Remark \ref{oss:orderonmu} $\ve{\rho}$ is strictly triangular, i.e., $\rho^k$ depends only on $\vxu$ and $\xv^1, \ldots, \xv^{k-1}$.
Hence, considering the coordinate $\xv^k$ of the conjugacy relation $\ve{\Phi} \circ \ve{f}=\ve{\widetilde{f}} \circ \ve{\Phi}$, we get
\begin{align*}
\xv^k \circ \ve{\Phi} \circ \ve{f} (\vx) 
&= \mu^k \xv^k + \rho^k(\vxu, \xv^1, \ldots, \xv^{k-1}) + R^k(\vx) + \phi^k \circ \ve{f} (\vx) \\ 
\xv^k \circ \ve{\widetilde{f}} \circ \ve{\Phi} (\vxu, \vxy, \vxv, \vxz) 
&= \mu^k \xv^k + \mu^k \phi^k(\vx) + \rho^k(\vxu, \xv^1, \ldots, \xv^{k-1})\mbox{.} 
\end{align*}
So we have to solve
$$
R^k(\vx) + \phi^k \circ \ve{f}(\vx)=\mu^k \phi^k(\vx)\mbox.
$$
It has an explicit solution, given by
$$
\phi^k(\vx)= \sum_{n=1}^\infty (\mu^k)^{-n} R^k \circ \ve{f}^{\comp n-1}(\vx)\mbox.
$$

Notice that for $\norm{\vx}$ small enough we have $\norm{\ve{f}^{\comp n}(\vx)} \leq \Lambda^n \norm{\vx}$. Then we have
$$
\abs{\phi^k(\vx)} \leq \sum_{n=1}^\infty \abs{\mu^k}^{-n} \abs{R^k \circ f^{\comp n-1}(\vx)} \leq \sum_{n=0}^\infty M \Lambda^{Nn} \abs{\mu^k}^{-n-1} \norm{\vx}\mbox,
$$
that converges since $\frac{\Lambda^N}{\abs{\mu^k}} < 1$.
\end{mystep}
\end{mydim}

\begin{myoss}\label{oss:nonarchiprimary}
The arguments of the proof of Theorem \ref{thm:rigidprimaryresonance} are also valid over any complete metrized field $\nK$.
Indeed, since \eqref{eqn:primarylot} is a linear (affine) equation on $\cc{\phi}{k}{\vi{n}}$ and $\cc{\widetilde{\rho}}{k}{(\vi{n}_\vxu, \vi{n}_\vxv)}$, it can be solved as well if $\nK$ is not algebraically closed.
Moreover, the estimates in Step $2$ works as well (or even better) in the non-archimedean case as in the complex case.

So Theorem \ref{thm:rigidprimaryresonance} holds in general, provided that all eigenvalues of $df_0$ belong to $\nK$ (see Remark \ref{oss:nonarchijordan}).
\end{myoss}

\mysect{Secondary Resonances}
\mysubsect{Resonance Relation}

\begin{mytext}
Starting from a germ written as in \eqref{eqn:rigidsplitprimaryresonance}, we can define $\vxx=(\vxu, \vxv)$, so that a contracting rigid germ $f:(\nC^d,0)\rightarrow (\nC^d, 0)$ is holomorphically conjugated to a map of the form
\begin{equation}\label{eqn:rigidprimaryresonance}
(\vxx, \vxy, \vxz) \mapsto \Big(\ve{\gamma} \vxx^\mi{P}+\ve{\sigma}(\vxx), \ve{\beta} \vxx^\mi{E}\vxy^\mi{D} \big(\one +  \ve{g}(\vxx,\vxy,\vxz)\big), \ve{h}(\vxx,\vxy,\vxz) \Big)\mbox,
\end{equation}
where
\begin{itemize}
\item $\vxx \in \nC^s$, $\vxy \in \nC^p$, and $\vxz \in \nC^{d-(s+p)}$;
\item $\ve{\gamma} \in (\nC^*)^s$ and $\ve{\beta} \in (\nC^*)^p$;
\item $\mi{P} \in \nmat{s}{s}{\nN}$ is a permutation matrix, $\mi{E} \in \nmat{s}{p}{\nN}$ and $\mi{D} \in \nmat{p}{p}{\nN}$;
\item $\ve{\sigma}:(\nC^s,0) \rightarrow (\nC^s, 0)$, $\ve{g}:(\nC^d,0) \rightarrow (\nC^p, 0)$ and $\ve{h}:(\nC^d,0) \rightarrow (\nC^{d-(s+p)}, 0)$;
\item $\{\vxy^\one=0\} \subseteq \critinf{f} \subseteq \{\vxx^\one \vxy^\one=0\}$;
\item $\ve{h}|_{\{\vxx=\vxy=\ve{0}\}}$ has nilpotent linear part;
\item $\ve{\sigma}$ is a polynomial map with only primary resonant monomials.
\end{itemize}
\end{mytext}
\begin{myoss}\label{oss:xuv}
The relation between equations \eqref{eqn:rigidsplitprimaryresonance} and \eqref{eqn:rigidprimaryresonance} is given by the identities:
$$
\mi{P}=
\begin{mymatrix}{2}
\mi{B} & \mi{0} \\
\mi{0} & \mi{\on{Id}}_e
\end{mymatrix}
\mbox,
\quad
\mi{E}=
\begin{mymatrix}{2}
\mi{C} & \mi{0} 
\end{mymatrix}
\mbox,
\quad
\ve{\gamma}=(\ve{\alpha},\ve{\mu})
\quad
\mbox{and }
\ve{\sigma}=(0, \ve{\rho})\mbox.
$$
\end{myoss}
\begin{mytext}
The aim of this section is to kill as many coefficients of $\ve{g}$ (expanded in formal power series) as possible, under the assumption of $\det \mi{D} \neq 0$ (i.e., injective internal action).
New formal obstructions appear: \myemph{secondary resonances}.
\end{mytext}

\begin{mydef}\label{def:secondaryresonances}
Let $f:(\nC^d,0) \rightarrow (\nC^d,0)$ be a contracting rigid germ as in \eqref{eqn:rigidprimaryresonance} with injective internal action, and let $\eta \in \nN^*$ be the order of $\mi{P}$.
A monomial $\vxx^{\vi{n}}$ is called \mydefin{secondary resonant} if
\begin{equation}\label{eqn:defsecondaryresonance}
\ve{\lambda}^{\eta\vi{n}} \in \on{Spec}(\mi{D}^\eta)\mbox,
\end{equation}
where $\ve{\lambda} \in (\nD^*)^s$ is the vector of non-zero eigenvalues of $df_0$ (counted with multiplicities).
\end{mydef}

\begin{myoss}
If $f:(\nC^d,0) \rightarrow (\nC^d,0)$ is a contracting rigid germ as in \eqref{eqn:rigidprimaryresonance}, and all the periodic irreducible components of $\critinf{f}$ are fixed, then $\eta=1$ and this definition coincides with the resonance relation given in the introduction.
\end{myoss}

\begin{mylem}
Let $\ve{f}:(\nC^d, 0) \rightarrow (\nC^d,0)$ be a contracting rigid germ written as in \eqref{eqn:rigidprimaryresonance}, with injective internal action.
Then there are only finitely-many secondary resonant monomials.
\end{mylem}
\begin{mydim}
It follows since $\mi{D}^\eta$ has only a finite number of eigenvalues $\mu$, and the secondary resonance relation is perfectly analogous to the primary resonance relation \eqref{eqn:defprimaryresonance}.
\end{mydim}

\begin{myes}
Let us see an example of how to compute secondary resonances.
Let $f:(\nC^3, 0) \rightarrow (\nC^3, 0)$ be a contracting rigid germ, with internal action $A$ given by
$$
\mi{A} =
\left(
\begin{array}{c|cc}
1 & 1 & 2 \\
\hline
0 & 2 & 1 \\
0 & 1 & 0
\end{array}
\right)
\mbox,
$$
where the splitting is according to the notations in \eqref{eqn:definternalaction}. Here
$$
\mi{D}=
\begin{mymatrix}{2}
2&1\\
1&0
\end{mymatrix}
\mbox,
$$
whose eigenvalues are $1 \pm \sqrt{2}$.

Set $\lambda$ the non-zero eigenvalue for $df_0$ and $(\xx^1, \xy^1, \xy^2)$ suitable coordinates in $0\in\nC^3$.
Then in this case $(\xx^1)^n$ is secondary resonant if
$$
\lambda^n =1-\sqrt{2}\mbox.
$$
Notice that $1+\sqrt{2}>1$ gives no resonances, being $\abs{\lambda}<1$.
\end{myes}

\begin{myoss}
We notice that secondary resonances for a contracting rigid germ $f:(\nC^d, 0) \rightarrow (\nC^d, 0)$ can appear only for $d \geq 3$; secondary resonances with periodic non-fixed irreducible components for $\critinf{f}$, or equivalently for $\eta \geq 2$ in \eqref{eqn:defsecondaryresonance}, can appear only for $d \geq 4$.
Primary and secondary resonances can appear in the same germ only for $d \geq 4$, and with $\eta \geq 2$ for $d \geq 5$.
\end{myoss}

\mysubsect{Main Theorem}

\begin{mytext}
Here we prove that we can kill all coefficients of $\ve{g}$ in \eqref{eqn:rigidprimaryresonance} except for secondary resonant monomials.
This theorem is the generalization of Theorem \ref{thm:rigidresonantpart} stated in the introduction.
\end{mytext}

\begin{mythm}\label{thm:rigidsecondaryresonance}
Let $\ve{f}:(\nC^d, 0) \rightarrow (\nC^d,0)$ be a contracting rigid germ with injective internal action.
Then $f$ is \myemph{analytically} conjugated to
\begin{equation}\label{eqn:rigidsecondaryresonance}
(\vxx, \vxy, \vxz) \mapsto \Big(\ve{\gamma} \vxx^\mi{P}+\ve{\sigma}(\vxx), \ve{\beta} \vxx^\mi{E}\vxy^\mi{D} \big(\one +  \ve{g}(\vxx)\big), \ve{h}(\vxx,\vxy,\vxz) \Big)\mbox,
\end{equation}
where
\begin{itemize}
\item $\vxx \in \nC^s$, $\vxy \in \nC^p$, and $\vxz \in \nC^{d-(s+p)}$;
\item $\ve{\gamma} \in (\nC^*)^s$ and $\ve{\beta} \in (\nC^*)^p$;
\item $\mi{P} \in \nmat{s}{s}{\nN}$ is a permutation matrix, $\mi{E} \in \nmat{s}{p}{\nN}$ and $\mi{D} \in \nmat{p}{p}{\nN}$ with $\det \mi{D} \neq 0$;
\item $\ve{\sigma}:(\nC^s,0) \rightarrow (\nC^s, 0)$, $\ve{g}:(\nC^s,0) \rightarrow (\nC^p, 0)$ and $\ve{h}:(\nC^d,0) \rightarrow (\nC^{d-(s+p)}, 0)$;
\item $\{\vxy^\one=0\} \subseteq \critinf{f} \subseteq \{\vxx^\one \vxy^\one=0\}$;
\item $\ve{h}|_{\{\vxx=\vxy=\ve{0}\}}$ has nilpotent linear part;
\item $\ve{\sigma}$ is a polynomial map with only primary resonant monomials,
\item $\ve{g}$ is a polynomial map with only secondary resonant monomials.
\end{itemize}
\end{mythm}
\begin{mydim}
We first prove in Step $1$ the formal counterpart of this theorem, and then we will deal with the convergence of the formal power series involved in Step $2$.
\begin{mystep}[Step $1$]
First of all, we can suppose that $f$ is of the form \eqref{eqn:rigidprimaryresonance}, with $\ve{h}|_{\{\vxx=\vxy=\ve{0}\}}$ that has a nilpotent lower triangular linear part.

We want to conjugate $f$ with a map $\widetilde{f}:(\nC^d, 0)\rightarrow (\nC^d, 0)$ of the form \eqref{eqn:rigidsecondaryresonance} (with $\widetilde{g}$ and $\widetilde{h}$ instead of $g$ and $h$ respectively).

Let us consider a local diffeomorphism $\Phi:(\nC^d, 0) \rightarrow (\nC^d, 0)$ of the form
$$
\ve{\Phi}(\vxx, \vxy, \vxz)=\Big(\vxx, \vxy \big(\one +  \ve{\phi}(\vxx, \vxy, \vxz)\big), \vxz\Big)\mbox,
$$
with $\phi:(\nC^d, 0)\rightarrow (\nC^p,0)$ a formal map.

Considering the conjugacy relation $\Phi \circ f = \widetilde{f} \circ \Phi$ for the coordinate $\vxy$, we get
\begin{align*}
\vxy \circ \ve{\Phi} \circ \ve{f} (\vxx, \vxy, \vxz) &=
\ve{\beta} \vxx^\mi{E}\vxy^\mi{D}\big(\one +  \ve{g}(\vxx, \vxy, \vxz)\big)\big(\one +  \ve{\phi} \circ \ve{f}(\vxx, \vxy, \vxz)\big)\\
\vxy \circ \ve{\widetilde{f}} \circ \ve{\Phi} (\vxx, \vxy, \vxz) &= 
\ve{\beta} \vxx^\mi{E}\vxy^\mi{D}\big(\one +  \ve{\phi}(\vxx, \vxy, \vxz)\big)^\mi{D} \big(\one +  \ve{\widetilde{g}} (\vxx)\big)
\mbox{.}
\end{align*}
Hence we have to solve
\begin{equation}\label{eqn:secondaryeqntosolve}
\big(\one +  \ve{g}(\vxx, \vxy, \vxz)\big)\big(\one +  \ve{\phi} \circ \ve{f}(\vxx, \vxy, \vxz)\big) = \big(\one +  \ve{\phi}(\vxx, \vxy, \vxz)\big)^\mi{D} \big(\one +  \ve{\widetilde{g}} (\vxx)\big)\mbox.
\end{equation}
Let us denote by $\lhs$ and $\rhs$ the left and right hand side of \eqref{eqn:secondaryeqntosolve} respectively.

We want now to expand in formal power series \eqref{eqn:secondaryeqntosolve} and to solve it defining (inductively) the coefficients of $\phi$ and $\widetilde{g}$.
Set $\vx=(\vxx, \vxy, \vxz)$ and
\begin{itemize}
\item $\vxx=(\xx^1, \ldots, \xx^s)$;
\item $\ve{\sigma}=(\sigma^1, \ldots, \sigma^s)$ and $\gamma^k \xx^k + \sigma^k(\vxx)=\sum^{\vi{n}} \cc{\sigma}{k}{\vi{n}_\vxx}\vxx^{\vi{n}_\vxx}$ for $1 \leq k \leq s$;
\item $\ve{\phi}=(\phi^1, \ldots, \phi^p)$ and $1 + \phi^k(\vx)= \sum_{\vi{n}} \cc{\phi}{k}{\vi{n}}\vx^\vi{n}$ for $1 \leq k \leq p$;
\item $\lhs=\big(\lhs^1, \ldots, \lhs^p\big)$ and $\lhs^k(\vx)= \sum_{\vi{n}} \cc{\lhs}{k}{\vi{n}}\vx^\vi{n}$ for $1 \leq k \leq p$, and analogously for $\rhs$;
\item $\ve{g}=(g^1, \ldots, g^p)$ and $1 + g^k(\vx)= \sum_{\vi{n}} \cc{g}{k}{\vi{n}}\vx^\vi{n}$ for $1 \leq k \leq p$;
\item $\ve{\widetilde{g}}=(\widetilde{g}^1, \ldots, \widetilde{g}^p)$ and $1 + \widetilde{g}^k(\vxx)= \sum_{\vi{n}_\vxx} \cc{\widetilde{g}}{k}{\vi{n}_\vxx}\vxx^{\vi{n}_\vxx}$ for $1 \leq k \leq p$;
\item $\ve{h}=(h^1, \ldots, h^{d-(s+p)})$ and $h^k(\vx)= \sum_{\vi{n}} \cc{h}{k}{\vi{n}}\vx^\vi{n}$ for $1 \leq k \leq d-(s+p)$.
\end{itemize}

Again, we split multi-indices $\vi{n}=(\vi{n}_\vxx, \vi{n}_\vxy, \vi{n}_\vxz) \in \nN^d$, where $\ _\vxx$ is the projection onto the coordinate $\vxx$, and similarly for other coordinates. In particular $\vi{n}_\vxx \in \nN^s$, $\vi{n}_\vxy \in \nN^p$ and $\vi{n}_\vxz \in \nN^{d-(s+p)}$.

By direct computations (see Remarks \ref{oss:expproductrule} and \ref{oss:powerseriesnotations}), we get
\begin{align}
\lhs^k&=\big(1+g^k(\vx)\big) \sum_{\vi{i} \in \nN^d} \left[\cc{\phi}{k}{\vi{i}}
\left(\ve{\gamma}\vxx^\mi{P} + \ve{\sigma}(\vxx)\right)^{\vi{i}_\vxx}
\left(\ve{\beta}\vxx^\mi{E}\vxy^\mi{D} \big(\one +  \ve{g}(\vx)\big)\right)^{\vi{i}_\vxy}
\left(\ve{h}(\vx)\right)^{\vi{i}_\vxz}
\right]\nonumber
\\
&= \sum_{\vi{i} \in \nN^d} \cc{\phi}{k}{\vi{i}}
\ve{\beta}^{\vi{i}_\vxy} \vxx^{\mi{E} \vi{i}_\vxy} \vxy^{\mi{D} \vi{i}_\vxy}
\left(
\sum_{\vpp{I} \in \nindd{\vi{i}_\vxx}{s}} \sigma_\vpp{I} \vxx^{\abs{\vpp{I}}}
\sum_{\vpp{J} \in \nind{\vi{i}_\vxy + \vi{e}^k}} g_\vpp{J} \vx^{\abs{\vpp{J}}}
\sum_{\vpp{K} \in \nind{\vi{i}_\vxz}} h_\vpp{K} \vx^{\abs{\vpp{K}}}
\right)\label{eqn:rsrlhsexpr}
\\ \nonumber \\
\rhs^k&=
\big(\one +  \ve{\phi}(\vx)\big)^{\mi{D} \vi{e}^k}
\sum_{\vi{j} \in \nN^s} \cc{\widetilde{g}}{k}{\vi{j}} \vxx^{\vi{j}} 
= \sum_{\vpp{H} \in \nind{\mi{D} \vi{e}^k}} \phi_\vpp{H} \vx^{\abs{\vpp{H}}}
\sum_{\vi{j} \in \nN^s} \cc{\widetilde{g}}{k}{\vi{j}} \vxx^{\vi{j}} \mbox, \label{eqn:rsrrhsexpr}
\end{align}
for $k=1, \ldots, p$, where $\vi{e}^k$ denotes the vector in $\nN^p$ with $1$ in the $k$-th position, and $0$ elsewhere.

Expressing explicitly the coefficients of $\lhs^k$ and $\rhs^k$ expanded in formal power series, from \eqref{eqn:rsrlhsexpr} and \eqref{eqn:rsrrhsexpr} respectively we obtain:
$$
\cc{\lhs}{k}{\vi{n}} = 
\sumaa{\substack{\vi{i}\in \nN^d \\ \vpp{I} \in \nindd{\vi{i}_\vxx}{s}, \vpp{J} \in \nind{\vi{i}_\vxy + \vi{e}^k}, \vpp{K} \in \nind{\vi{i}_\vxz}\\ \cond{1}}}{1.2cm} 
\cc{\phi}{k}{\vi{i}} \ve{\beta}^{\vi{i}_\vxy} \sigma_{\vpp{I}} g_\vpp{J}  h_{\vpp{K}}
\mbox, 
\qquad
\cc{\rhs}{k}{\vi{n}} =
\sumaa{\substack{\vi{j}\in \nN^s \\ \vpp{H} \in \nind{\mi{D} \vi{e}^k}\\ \cond{2}}}{0.2cm}
\cc{\widetilde{g}}{k}{\vi{j}}\phi_{\vpp{H}}
\mbox,
$$
where
$$
\cond{1}=
\left\{
\begin{array}{l}
\mi{E} \vi{i}_\vxy +\abs{\vpp{I}} +\abs{\vpp{J}}_\vxx +\abs{\vpp{K}}_\vxx= \vi{n}_\vxx \\
\mi{D} \vi{i}_\vxy +\abs{\vpp{J}}_\vxy +\abs{\vpp{K}}_\vxy= \vi{n}_\vxy \\
\abs{\vpp{J}}_\vxz +\abs{\vpp{K}}_\vxz= \vi{n}_\vxz
\end{array}
\right.
\mbox,
$$
and
$$
\cond{2}=
\left\{
\begin{array}{l}
\vi{j} + \abs{\vpp{H}}_\vxx = \vi{n}_\vxx \\
\abs{\vpp{H}}_\vxy = \vi{n}_\vxy \\
\abs{\vpp{H}}_\vxz = \vi{n}_\vxz 
\end{array}
\right.
\mbox.
$$

We want to solve the equation 
\begin{equation}\label{eqn:secondarycoefftosolve}
\cc{\eq}{k}{\vi{n}}:= \cc{\rhs}{k}{\vi{n}}-\cc{\lhs}{k}{\vi{n}}=0
\end{equation}
for every $k$ and $\vi{n}$, with respect to the coefficients $\cc{\phi}{k}{\vi{n}}$ of $\ve{\phi}$ and $\cc{\widetilde{g}}{k}{\vi{n}_\vxx}$ of $\ve{\widetilde{g}}$.

We recall the partial order and the total order on indices in $\nN^d$ that we need to make computations.
Set $\vi{n}=(n_1, \ldots, n_d)$ and $\vi{m}=(m_1, \ldots, m_d)$.
\begin{trivlist}
\item[Partial order $\preceq$:] we say that $\vi{m} \preceq \vi{n}$ iff we have $m_k \leq n_k$ for every $k=1, \ldots, d$.
\item[Total order $\leq$:] we say that $\vi{m} \leq \vi{n}$ iff $(\abs{\vi{m}}, m_1, \ldots, m_d) \leq_{\on{lex}} (\abs{\vi{n}}, n_1, \ldots, n_d)$, where $\leq_{\on{lex}}$ is the lexicographic order (on $\nN^{d+1}$).
\end{trivlist}

\begin{mydef}
Let $k \in \{1, \ldots, p\}$ be an integer and $\vi{n}$ be a multi-index (in $\nN^d$ or $\nN^s$).
We call \mydefin{weight} of $(k, \vi{n})$ the value
$$
\wa{k}{\vi{n}}:=\abs{\vi{n}} \in \nN\mbox.
$$
As in Step $1$ of the proof of Theorem \ref{thm:rigidprimaryresonance}, the notation 
$$
\llot{\abs{\vi{n}}}{\phi, \widetilde{g}}
$$
stands for a suitable polynomial in $\cc{\phi}{l}{\vi{m}}$ and $\cc{\widetilde{g}}{l}{\vi{m}_\vxx}$ satisfying 
$$
\wa{l}{\vi{m}} < \wb{\abs{\vi{n}}}\mbox.
$$
We shall also omit $\widetilde{g}$ when the polynomial does not depend on any coefficient $\cc{\widetilde{g}}{l}{\vi{m}_\vxx}$.
\end{mydef}

Notice that the definition of weight here is slightly different from the one given by Definition \ref{def:primaryweight}.
Still, we have that for every $W \in \nN$, there are only finitely many $(k, \vi{n})$ such that $\wa{k}{\vi{n}} \leq W$.

\begin{mylem}\label{lem:secondaryrhslot}
For every $k=1,\ldots, p$ and $\vi{n} \neq \vi{0}$ we have
\begin{equation}\label{eqn:secondaryrhslot}
\cc{\rhs}{k}{\vi{n}} = \kronecker{\vi{n}_\vxy}{\vi{0}} \kronecker{\vi{n}_\vxz}{\vi{0}} \cc{\widetilde{g}}{k}{\vi{n}_\vxx} +
\sum_{l=1}^p d_l^k \cc{\phi}{l}{\vi{n}}
+ \llot{\abs{\vi{n}}}{\phi, \widetilde{g}}\mbox,
\end{equation}
where $\kronecker{}{}$ denotes the Kronecker's delta function and $\mi{D}=(d_l^k)$.
\end{mylem}
\begin{mydim}
Set $W=\wb{\abs{\vi{n}}}$.
From the first equation of $\cond{2}$ we have $\vi{j} \preceq \vi{n}_\vxx$.
Hence the only term of the form $\cc{\widetilde{g}}{k}{\vi{j}}$ whose weight is $\geq W$ is given by $\vi{j}=\vi{n}_\vxx$, when $\vi{n}_\vxy=\vi{0}$ and $\vi{n}_\vxz=\vi{0}$.
In this case, $\abs{\vpp{H}}=\vi{0}$, and $\phi_\vpp{H}=1$ (being $\cc{\phi}{k}{\vi{0}}=1$ for $k=1, \ldots p$).
This gives the first term of \eqref{eqn:secondaryrhslot}.

Since $\abs{\vpp{H}} \preceq \vi{n}$, the only terms $\cc{\phi}{l}{\vi{m}}$ with $\wa{l}{\vi{m}}\geq W$ that appear are when $\vi{m}=\vi{n}$, and
$$
\vpp{H}=
\big(
\underbrace{\vi{0}, \ldots, \vi{0}}_{d_1^k} |
\cdots |
\underbrace{\vi{0}, \ldots, \vi{0}, \vi{n}, \vi{0}, \ldots, \vi{0}}_{d_l^k} |
\cdots |
\underbrace{\vi{0}, \ldots, \vi{0}}_{d_p^k}
\big)
\mbox.
$$
Since we have $d_l^k$ choices for where to put $\vi{n}$, \eqref{eqn:secondaryrhslot} follows.
\end{mydim}

\begin{mylem}\label{lem:secondarylhslot}
For every $k=1,\ldots, p$ and $\vi{n} \neq \vi{0}$ we have
\begin{equation}\label{eqn:secondarylhslot}
\cc{\lhs}{k}{\vi{n}} = \kronecker{\vi{n}_\vxy}{\vi{0}}\kronecker{\vi{n}_\vxz}{\vi{0}} \ve{\gamma}^{\mi{P}^{-1} \vi{n}_\vxx} \cc{\phi}{k}{(\mi{P}^{-1} \vi{n}_\vxx, \vi{0}, \vi{0})}  + \llot{\abs{\vi{n}}}{\phi}\mbox.
\end{equation}
\end{mylem}
\begin{mydim}
Thanks to Lemma \ref{lem:computationpsiH}, we get that $\sigma_{\vpp{I}}\neq 0$ only if $\abs{\vpp{I}} \geq \vi{i}_\vxx$.
Lemma \ref{lem:computationpsiH}.(i) says that $h_{\vpp{K}}\neq 0$ only if $\abs{\vpp{K}} > (\vi{0}, \vi{0}, \vi{i}_\vxz)$ when $\vi{i}_\vxz \neq \vi{0}$.
Moreover, we have $\abs{\vi{i}_\vxy} < \abs{\mi{D} \vi{i}_\vxy}$ if $\vi{i}_\vxy \neq 0$.
Then we have
$$
\abs{\vi{i}} = \abs{\vi{i}_\vxx} + \abs{\vi{i}_\vxy} + \abs{\vi{i}_\vxz} \leq \abs{\abs{\vpp{I}}} + \abs{\vi{i}_\vxy} + \abs{\abs{\vpp{K}}} = \abs{\vi{n}} - \abs{\abs{\vpp{J}}} - \abs{\mi{E} \vi{i}_\vxy} - \abs{\mi{D} \vi{i}_\vxy} + \abs{\vi{i}_\vxy} \leq \abs{\vi{n}}\mbox,
$$
where the equality can hold only when $\vi{i}_\vxy = \vi{0}$ and $\abs{\abs{\vpp{J}}}=0$.
Suppose this is the case; then $\vpp{J}$ is made by just one column (in position $k$) made by $0$'s, and hence $g_\vpp{J}=\cc{g}{k}{\vi{0}}=1$.
From the first equation of $\cond{1}$ we also get that $\vi{i}_\vxx \leq \abs{\vpp{I}} \preceq \vi{n}_\vxx$. It follows that the only terms $\cc{\phi}{k}{\vi{i}}$ whose weight is $\geq \abs{\vi{n}}$ appear when $\abs{\vi{i}_\vxx}=\abs{\vi{n}_\vxx}$, $\abs{\vpp{I}}=\vi{n}_\vxx$ and $\vi{n}_\vxy=\vi{0}$.

In this case ($\vi{i}_\vxy=\vi{n}_\vxy=\vi{0}$, $\vpp{J}=\vi{0} \in \nind{\vi{e}^k}$ and $\abs{\vpp{I}}=\vi{n}_\vxx$), $\cond{1}$ becomes
$$
\left\{
\begin{array}{l}
\abs{\vpp{K}}_\vxx= \vi{0} \\
\abs{\vpp{K}}_\vxy= \vi{0} \\
\abs{\vpp{K}}_\vxz= \vi{n}_\vxz
\end{array}
\right.
\mbox.
$$
From Lemma \ref{lem:computationpsiH}.(i), being $\abs{K}=(\vi{0}, \vi{0}, \abs{\vpp{K}}_\vxz)$, it follows that the only term with weight $\geq \wb{\abs{n}}$ appear when in addition $\vi{i}_\vxz=\vi{n}_\vxz=\vi{0}$.
In this case $\vpp{K}=\emptyset$ and $h_\emptyset=1$.

We shall show now that the conditions $\abs{\vpp{I}} = \vi{n}_\vxx$ and $\sigma_\vpp{I} \neq 0$ are satisfied by a unique $\vpp{I} \in \nindd{\vi{i}_\vxx}{s}$, and in this case $\sigma_\vpp{I}=\gamma^{\vi{P}^{-1} \vi{n}_\vxx}$.

Let us split again $\vxx=(\vxu, \vxv)$, $\ve{\sigma}=(\ve{0}, \ve{\rho})$, $\mi{P}=\on{Diag}(\mi{B}, \idmat)$ and $\gamma=(\alpha, \mu)$ as in Remark \ref{oss:xuv}.
The condition $\abs{\vpp{I}}=\vi{n}_\vxx$ becomes
$$
\left\{
\begin{array}{l}
\mi{B} \vi{i}_\vxu = \vi{n}_\vxu \\
\abs{\vpp{I}}_\vxv = \vi{n}_\vxv
\end{array}
\right.
$$
and $\sigma_\vpp{I}=\alpha^{\vi{i}_\vxu} \rho_{\vpp{I}_\vxv}$.
Then $\vi{i}_\vxu=\vi{B}^{-1} \vi{n}_\vxu$, and thanks to Lemma \ref{lem:computationpsiH}.(ii) we get $\rho_{\vpp{I}_\vxv}=\mu^{\vi{n}_\vxv}$.

Writing again with the previous notations, we get the statement.
\end{mydim}

Set $\eq_\vi{n}:=(\cc{\eq}{1}{\vi{n}},\ldots, \cc{\eq}{p}{\vi{n}})$,
$\ve{\phi}_{\vi{n}}:=(\cc{\phi}{1}{\vi{n}}, \ldots, \cc{\phi}{p}{\vi{n}})$ and
$\ve{\widetilde{g}}_{\vi{n}_\vxx}:=(\cc{\widetilde{g}}{1}{\vi{n}_\vxx}, \ldots, \cc{\widetilde{g}}{p}{\vi{n}_\vxx})$.
Thanks to Lemmas \ref{lem:secondaryrhslot} and \ref{lem:secondarylhslot}, $\eq_\vi{n}=\ve{0}$ becomes
\begin{equation}\label{eqn:secondarylot}
\ve{\phi}_{\vi{n}} \mi{D} + \kronecker{\vi{n}_\vxy}{\vi{0}}\kronecker{\vi{n}_\vxz}{\vi{0}} \big(\ve{\widetilde{g}}_{\vi{n}_\vxx} - \ve{\gamma}^{\vi{n}_\vxx} \ve{\phi}_{\vi{n}}\big) = \llot{\abs{\vi{n}}}{\phi, \widetilde{g}}\mbox.
\end{equation}
This affine equation, where the unknowns are $\ve{\phi}_{\vi{n}}$ and $\ve{\widetilde{g}}_{\vi{n}_\vxx}$, has always a solution.
At this point we conjugated $f$ to a map $\widetilde{f}$ as in \eqref{eqn:rigidsecondaryresonance}, but with $\widetilde{g}:(\nC^s,0) \rightarrow (\nC^p,0)$ a (vector of) formal power series.
Next we show that we can solve \eqref{eqn:secondarylot} and get $\widetilde{g}$ polynomial with only secondary resonant monomials.

We solve $\eq_\vi{n}=\ve{0}$ inductively on $\wb{\abs{\vi{n}}}$ as follows.
\begin{trivlist}
\item If $\wb{\abs{\vi{n}}} = 0$, i.e., if $\vi{n}=\ve{0}$, we set $\ve{\phi}_{\vi{0}}:=\one$ and $\ve{\widetilde{g}}_{\vi{0}_\vxx}=\one$.
\item Set $0<W \in \nN$, and suppose that $\ve{\phi}_{\vi{m}}$ and $\ve{\widetilde{g}}_{\vi{m}_\vxx}$ are known for $\abs{\vi{m}} < W$.
We want to solve \eqref{eqn:secondarylot} for $\abs{\vi{n}}=W$.

Notice that $\llot{\abs{\vi{n}}}{\phi, \widetilde{g}}$ is a polynomial that depend on $\cc{\phi}{l}{\vi{m}}$ and $\cc{\widetilde{g}}{l}{\vi{m}_\vxx}$ only for weights strictly less than $W$.
Hence thanks to the induction hypothesis, $\llot{\abs{\vi{n}}}{\phi, \widetilde{g}}$ is a known value in $\nC^p$.
\begin{trivlist}
\item Suppose $(\vi{n}_\vxy, \vi{n}_\vxz) \neq (\vi{0}, \vi{0})$. Then \eqref{eqn:secondarylot} becomes
$$
\ve{\phi}_{\vi{n}} \mi{D} = \llot{\abs{\vi{n}}}{\phi, \widetilde{g}}\mbox.
$$
and being $\det \mi{D} \neq 0$, there exists a unique $\phi_\vi{n} \in \nC^p$ that solves the equation.
\item Suppose $(\vi{n}_\vxy, \vi{n}_\vxz) = (\vi{0}, \vi{0})$. Then \eqref{eqn:secondarylot} becomes
\begin{equation}\label{eqn:secondarylot00}
- \ve{\phi}_{\vi{n}} \mi{D} + \ve{\gamma}^{\mi{P}^{-1}\vi{n}_\vxx} \ve{\phi}_{(\mi{P}^{-1}\vi{n}_\vxx, \vi{0}, \vi{0})} = \ve{\widetilde{g}}_{\vi{n}_\vxx} + \llot{\abs{\vi{n}}}{\phi, \widetilde{g}}\mbox.
\end{equation}
Suppose that $\vi{n}_\vxu=\mi{P}^{-1} \vi{n}_\vxu$: we have two cases.

Suppose $\mi{D} - \ve{\gamma}^{\vi{n}_\vxx} \idmat$ is invertible, i.e., $\vxx^{\vi{n}_\vxx}$ is not secondary resonant.
Then we can put $\ve{\widetilde{g}}_{\vi{n}_\vxx}=0$ and there exists a unique $\ve{\phi}_{\vi{n}} \in \nC^p$ that solves the equation.

Suppose $\mi{D} - \ve{\gamma}^{\vi{n}_\vxx} \idmat$ is not invertible, i.e., $\vxx^{\vi{n}_\vxx}$ is secondary resonant.
Then we can put $\ve{\phi}_{\vi{n}}=0$ and there exists a unique $\ve{\widetilde{g}}_{\vi{n}_\vxx} \in \nC^p$ that solves the equation.

In the general case, let $\widetilde{\eta}$ be the smallest number in $\nN^*$ such that $\vi{n}_\vxx=\mi{P}^{\widetilde{\eta}} \vi{n}_\vxx$.
Set $\vi{n}_\vxx^{(l)}:=\mi{P}^{-l} \vi{n}_\vxx$ for $l=0, \ldots, \widetilde{\eta}-1$.
We consider the equation \eqref{eqn:secondarylot00} for $\vi{n}_\vxx, \vi{n}_\vxx^{(1)}, \ldots, \vi{n}_\vxx^{(\widetilde{\eta}-1)}$ simultaneously.

We get the following (block) linear system:
\begin{equation}\label{eqn:secondarylinearsystem}
\begin{mymatrix}{1}
\!\!\!\phi_{(\vi{n}_\vxx, \vi{0}, \vi{0})}\!\!\!\\
\!\!\!\phi_{(\vi{n}_\vxx^{(1)}, \vi{0}, \vi{0})}\!\!\!\\
\vdots \\
\vdots \\
\!\!\!\phi_{(\vi{n}_\vxx^{(\widetilde{\eta}-1)}, \vi{0}, \vi{0})}\!\!\!
\end{mymatrix}^T
\begin{mymatrix}{5}
-\mi{D} & 0 & \cdots & 0 & \Big.\!\!\!\ve{\gamma}^{\vi{n}_\vxx} \on{Id}_p \!\!\!\Big. \\
\!\!\!\ve{\gamma}^{\vi{n}_\vxx^{(1)}} \on{Id}_p \!\!\! & -\mi{D} & 0 & \ddots & 0 \\
0 & \!\!\! \ve{\gamma}^{\vi{n}_\vxx^{(2)}} \on{Id}_p \!\!\! & \ddots & \ddots & \vdots \\
\vdots & \ddots & \ddots & -\mi{D} & 0 \\
0 & \cdots & 0 & \!\!\! \ve{\gamma}^{\vi{n}_\vxx^{(\widetilde{\eta}-1)}} \on{Id}_p \!\!\! & -\mi{D}
\end{mymatrix}
=
\begin{mymatrix}{1}
\!\!\!\widetilde{g}_{\vi{n}_\vxx}\!\!\!\\
\!\!\!\widetilde{g}_{\vi{n}_\vxx^{(1)}}\!\!\!\\
\vdots \\
\vdots \\
\!\!\!\widetilde{g}_{\vi{n}_\vxx^{(\widetilde{\eta}-1)}}\!\!\!
\end{mymatrix}^T
+ \lot\mbox,
\end{equation}
where
$$
\lot=\llot{\abs{\vi{n}}}{\phi, \widetilde{g}}\mbox.
$$
for each coordinate.

Let us consider the linear combination of the columns (numbered from $1$ to $\widetilde{\eta}$) of the linear system \eqref{eqn:secondarylinearsystem}, where the $l$-th column is multiplied by
$$
\Delta_l:=\left(\prod_{h=1}^{l-1} \ve{\gamma}^{\vi{n}_\vxx^{(h)}}\right) \mi{D}^{\widetilde{\eta}-l}\mbox.
$$
Then we get
$$
\left(-\mi{D}^{\widetilde{\eta}} + \prod_{l=0}^{\widetilde{\eta}-1} \ve{\gamma}^{\vi{n}_\vxx^{(l)}} \on{Id}_p\right) \phi_{(\vi{n}_\vxx, \vi{0}, \vi{0})} = \sum_{l=1}^{\widetilde{\eta}} \Delta_l \widetilde{g}_{\vi{n}_\vxx^{(l-1)}} + \lot\mbox.
$$
Since $\det \mi{D} \neq 0$, it follows that the linear system \eqref{eqn:secondarylinearsystem} is invertible iff $\vi{n}_\vxx$ is not secondary resonant.

In this case we can put $\widetilde{g}_{\vi{n}_\vxx^{(l)}}=0$ for every $l=0, \ldots, \widetilde{\eta}-1$, and there exist (unique) $\phi_{(\vi{n}_\vxx^{(l)}, \vi{0}, \vi{0})} \in \nC^p$ for $l=0, \ldots, \widetilde{\eta}-1$ that satisfy \eqref{eqn:secondarylinearsystem}.

If $\vi{n}_\vxx$ is secondary resonant we can still set any value for $\phi_{(\vi{n}_\vxx^{(l)}, \vi{0}, \vi{0})}$ (for example, all equal to $0$), and find unique $\widetilde{g}_{\vi{n}_\vxx^{(l)}} \in \nC^p$ for every $l=0, \ldots, \widetilde{\eta}-1$ that satisfy \eqref{eqn:secondarylinearsystem}.
\end{trivlist}
\end{trivlist}

As in the proof of Theorem \ref{thm:rigidprimaryresonance}, we have defined the conjugation $\Phi$ as an invertible formal map so we can then define $\widetilde{h}$ such that the conjugacy relation \eqref{eqn:secondaryeqntosolve} holds for all coordinates.
 
\end{mystep}

\begin{mystep}[Step $2$]
The following estimations are quite standard.
Pick $0 < \Lambda < 1$ such that $\Lambda > \on{specrad}(df_0)$ the spectral radius of the differential $df_0$ of $f$ at $0$, and take $N$ big enough such that $\abs{\mi{D}^{-1}} \Lambda^N<1$ and no secondary resonances $\vxx^{\vi{n}}$ appear for $\abs{\vi{n}} \geq N$.

For proving the formal result, we introduced a weight, and noticed that for every $W \in \nN$, there are only finitely many $(k, \vi{n})$ such that $\wa{k}{\vi{n}} \leq W$.
It follows that there exist $M >0$ and a polynomial (hence holomorphic) change of coordinates that conjugates $\ve{f}$ with a map of the form
\begin{equation}\label{eqn:rigidsecondaryresonancewithrest}
(\vxx, \vxy, \vxz) \mapsto \Big(\ve{\gamma} \vxx^\mi{P}+\ve{\sigma}(\vxx), \ve{\beta} \vxx^\mi{E}\vxy^\mi{D} \big(\one +  \ve{g}(\vxx) + \ve{R}(\vxx,\vxy,\vxz)\big), \ve{h}(\vxx,\vxy,\vxz) \Big)\mbox,
\end{equation}
with the same conditions as for \eqref{eqn:rigidsecondaryresonance} and $R:(\nC^d, 0) \rightarrow (\nC^p, 0)$ such that
$$
\norm{\ve{R}(\vx)} \leq M \norm{\vx}^N
$$
for a suitable $M > 0$ and $\norm{\vx}$ small enough, where $\vx=(\vxx, \vxy, \vxz)$.

We can hence suppose that $f$ is of the form \eqref{eqn:rigidsecondaryresonancewithrest}, and try to kill the map $\ve{R}$: we look for a conjugacy between $f$ and a map $\widetilde{f}$ of the form \eqref{eqn:rigidsecondaryresonance} (with $\widetilde{h}$ instead of $h$).

Let us consider then a local diffeomorphism of the form
$$
\ve{\Phi}(\vxx, \vxy, \vxz)=\Big(\vxx, \vxy \big(\one +  \ve{\phi}(\vxx, \vxy, \vxz)\big), \vxz\Big)\mbox.
$$
Looking at the conjugacy relation $\Phi \circ f = \widetilde{f} \circ \Phi$ at the coordinate $\vxy$, we get
\begin{align*}
\vxy \circ \Phi \circ f (\vxx, \vxy, \vxz) 
&= \ve{\beta} \vxx^\mi{E} \vxy^\mi{D} \big(\one +  \ve{g}(\vxx) + \ve{R}(\vxx, \vxy, \vxz)\big)\big(\one +  \ve{\phi} \circ \ve{f}(\vxx, \vxy, \vxz)\big) \\ 
\vxy \circ \widetilde{f} \circ \Phi (\vxx, \vxy, \vxz) 
&= \ve{\beta} \vxx^\mi{E} \vxy^\mi{D} \big(\one +  \ve{\phi}(\vxx, \vxy, \vxz)\big)^\mi{D} \big(\one +  \ve{g}(\vxx)\big)\mbox{.} 
\end{align*}
Hence we have to solve
\begin{equation}\label{eqn:rigidsecondaryresonanceanalytictosolve}
\big(\one +  \ve{\phi}(\vx)\big)^\mi{D} = \big(\one +  \ve{\phi} \circ \ve{f}(\vx)\big) \big(\one +  \ve{e}(\vx)\big)
\mbox,
\end{equation}
where $\vx=(\vxx, \vxy, \vxz)$ and
$$
\ve{e}(\vxx, \vxy, \vxz)=\frac{R(\vxx, \vxy, \vxz)}{\one +  \ve{g}(\vxx)}\mbox.
$$
In particular we have
$$
\norm{\ve{e}(\vx)} \leq K \norm{\vx}^N
$$
for $K>0$ big enough and $\norm{\vx}$ small enough.

Equation \eqref{eqn:rigidsecondaryresonanceanalytictosolve} has an explicit solution, given by
\begin{equation*}
\one +  \ve{\phi}(\vx) = \prod_{n=1}^\infty \big(\one +  \ve{e} \circ \ve{f}^{\comp {n-1}}(\vx)\big)^{\mi{D}^{-n}}\mbox;
\end{equation*}
let us show that this product is convergent.

Thanks to Proposition \ref{prop:takingthelog}, we just need to prove that 
$$
\sum_{n=1}^\infty \left(\ve{e} \circ \ve{f}^{\comp n-1}(\vx)\right)\mi{D}^{-n}
$$
converges for $\norm{\vx}$ small enough.

Notice that for $\norm{\vx}$ small enough we have $\norm{\ve{f}^{\comp n}(\vx)} \leq \Lambda^n \norm{\vx}$.

Then we have
$$
\norm{\sum_{n=1}^\infty \left(\ve{e} \circ \ve{f}^{\comp n-1}(\vx)\right) \mi{D}^{-n}} \leq
\sum_{n=1}^\infty \abs{\mi{D}^{-n}} \norm{\ve{e} \circ \ve{f}^{\comp n-1}(\vx)} \leq
\sum_{n=1}^\infty \abs{\mi{D}^{-1}}^n K \Lambda^{(n-1)N} \norm{\vx}^N\mbox,
$$
that converges being $\abs{\mi{D}^{-1}} \Lambda^N < 1$.
\end{mystep}
\end{mydim}

\begin{myoss}\label{oss:noninjectiveinternalaction}
Let us take a rigid germ that has a non-injective internal action: we can write it in the form \eqref{eqn:rigidprimaryresonance}, with $\det \mi{D} = 0$ (suppose also $\mi{P}=\on{Id}$ for simplicity). We can try to kill, at least formally, as many coefficients of $\ve{g}$ as possible, as we did in the case of injective internal action.
Proceeding as in the proof of Theorem \ref{thm:rigidsecondaryresonance}, we get an equation to solve of the form \eqref{eqn:secondarylot}.
When $\vi{n}_\vxy$ or $\vi{n}_\vxz$ are different from $0$, the linear system becomes
\begin{equation*}
\ve{\phi}_{\vi{n}} \mi{D} = \llot{\abs{\vi{n}}}{\phi, \widetilde{g}}\mbox,
\end{equation*}
that is not invertible, being $\det\mi{D} \neq 0$.
So in general, besides the secondary resonances already described, some other resonances of the form $\vxx^{\vi{n}_\vxx}\vxy^{\vi{n}_\vxy}\vxz^{\vi{n}_\vxz}$ with $(\vi{n}_\vxy,\vi{n}_\vxz) \neq (0,0)$ will appear. 
\end{myoss}

\begin{myoss}\label{oss:nonarchisecondary}
Theorem \ref{thm:rigidsecondaryresonance} holds over any complete metrized field $\nK$ of characteristic $0$ (provided that all eigenvalues of $df_0$ belong to $\nK$).
The reasons are the same as for Theorem \ref{thm:rigidprimaryresonance} (see Remark \ref{oss:nonarchiprimary}).
The theorem fails, already for $d=p=1$, over a field of positive characteristic.
In fact, although $\mi{D}$ is invertible as a matrix with integer (rational) coefficients, it could not be invertible when seen as a matrix with coefficients in $\nK$.
If this is the case, equation \eqref{eqn:secondarylot} could not be solved in general.
\end{myoss}

\mysect{Rigid Germs with $s+p=d-1$}
Theorem \ref{thm:rigidsecondaryresonance} gives in particular the complete classification of contracting rigid germs with injective internal action such that $s+p=d$, where as before $s$ is the number of non-zero eigenvalues of $df_0$, and $p$ is the number of non-periodic components of $\critinf{f}$.

In this section we shall deal with the case of a contracting rigid germ with injective internal action such that $s+p=d-1$.
Thanks to Theorem \ref{thm:rigidsecondaryresonance} we can holomorphically conjugate $f$ with a map of the form \eqref{eqn:rigidsecondaryresonance}, with $h:(\nC^d, 0) \rightarrow (\nC, 0)$ and $\vxz \in \nC$.

In this case we can say more, and get a similar result of what happens in the $2$-dimensional case (see \cite[pp. 491--494]{favre:rigidgerms}). 

\begin{mythm}\label{thm:rigidaffine}
Let $f:(\nC^d, 0) \rightarrow (\nC^d,0)$ be a contracting rigid germ with injective internal action, and such that $s+p=d-1$, where $s$ is the number of non-zero eigenvalues of $df_0$, and $p$ is the number of non-periodic components of $\critinf{f}$.
Then $f$ is \myemph{analytically} conjugated to a map of the form
\begin{equation}\label{eqn:rigidaffinepart}
(\vxx, \vxy, \xz) \mapsto \Big(\ve{\gamma} \vxx^{\mi{P}} + \ve{\sigma}(\vxx), \ve{\beta} \vxx^\mi{E}\vxy^\mi{D} \big(\one + \ve{g}(\vxx)\big),\nu \vxx^\vi{l}\vxy^\vi{m} \xz + \omega(\vxx, \vxy)\Big)\mbox,
\end{equation}
where
\begin{itemize}
\item $\vxx \in \nC^s$, $\vxy \in \nC^p$ and $\xz \in \nC$;
\item $\ve{\gamma} \in (\nC^*)^s$, $\ve{\beta} \in (\nC^*)^p$ and $\nu \in \nC^*$;
\item $\mi{P} \in \nmat{s}{s}{\nN}$ is a permutation matrix, $\mi{E} \in \nmat{s}{p}{\nN}$, $\mi{D} \in \nmat{p}{p}{\nN}$ and $(\vi{l}, \vi{m}) \in \nN^s \times \nN^p \setminus \{(\ve{0},\ve{0})\}$;
\item $\ve{\sigma}:(\nC^s,0) \rightarrow (\nC^s, 0)$, $\ve{g}:(\nC^s,0) \rightarrow (\nC^p, 0)$ and $\omega:(\nC^{d-1},0) \rightarrow (\nC, 0)$;
\item $\{\vxy^\one=0\} \subseteq \critinf{f} \subseteq \{\vxx^\one \vxy^\one=0\}$;
\item $\ve{\sigma}$ is a polynomial map with only primary resonant monomials;
\item $\ve{g}$ is a polynomial map with only secondary resonant monomials;
\item $\omega$ is analytic.
\end{itemize}
For $d \geq 3$ we cannot get in general $\omega$ polynomial (see Remark \ref{oss:nonpolnormalform}).
\end{mythm}

\begin{myoss}\label{oss:rigidd-1reducible}
Let us suppose that $f:(\nC^d, 0) \rightarrow (\nC^d,0)$ is a contracting rigid germ as in \eqref{eqn:rigidsecondaryresonance} and satisfying the hypotheses of Theorem \ref{thm:rigidaffine}.

Let us split again $\vxx=(\vxu, \vxv)$ (see Remark \ref{oss:xuv}), with $\vxu \in \nC^r$ and $\vxv \in \nC^e$: then $f$ is of the form
\begin{equation}\label{eqn:rigidsecondaryresonancesplit}
(\vxu, \vxv, \vxy, \vxz) \mapsto \Big(\ve{\alpha} \vxu^\mi{B}, \ve{\mu} \vxv + \ve{\rho}(\vxu, \vxv), \ve{\beta} \vxu^\mi{C}\vxy^\mi{D} \big(\one + \ve{g}(\vxu,\vxv)\big), h(\vxu,\vxv,\vxy,\xz) \Big)\mbox.
\end{equation}
If we compute $\det df$, we get
$$
\det df = \vxu^{\vi{a}} \vxy^{\vi{b}} \frac{\partial h}{\partial \xz} U(\vx)\mbox,
$$
for suitable $\vi{a} \in \nN^r$, $\vi{b} \in \nN^p$ and a holomorphic map $U:\nC^d \rightarrow \nC$ with $U(\ve{0}) \neq 0$, where $\vx=(\vxu, \vxv, \vxy, \xz)$.

Since $\critinf{f}=\{\vxu^\one \vxy^\one=0\}$, we get
$$
\frac{\partial h}{\partial \xz}= \vxu^{\vi{l}_\vxu}\vxy^\vi{m} V(\vx)\mbox,
$$
with $\vi{l}_\vxu \in \nN^r$, $\vi{m} \in \nN^p$ and $V(\ve{0}) \neq 0$.
Integrating, we obtain
$$
h(\vxu, \vxv, \vxy, \xz)=\nu \vxu^{\vi{l}_\vxu}\vxy^\vi{m} \xz \big(1+\varepsilon(\vxu, \vxv, \vxy, \xz)\big) + \omega(\vxu, \vxv, \vxy)\mbox,
$$
with $\nu \in \nC^*$, $\varepsilon:(\nC^d, 0) \rightarrow (\nC, 0)$ and $\omega:(\nC^{d-1}, 0) \rightarrow (\nC, 0)$ (and $(\vi{l}_\vxu, \vi{m}) \neq \ve{0}$).

As in Remark \ref{oss:xuv}, to simplify notations we use $\vxx$ instead of $(\vxu, \vxv)$; summing up, we can suppose that $f$ is of the form
\begin{equation}\label{eqn:rigidd-1reducible}
(\vxx, \vxy, \xz) \mapsto \Big(\ve{\gamma} \vxx^{\mi{P}} + \ve{\sigma}(\vxx), \ve{\beta} \vxx^\mi{E}\vxy^\mi{D} \big(\one + \ve{g}(\vxx)\big),\nu \vxx^\vi{l}\vxy^\vi{m} \xz \big(1+\varepsilon(\vxx, \vxy, \xz)\big) + \omega(\vxx, \vxy)\Big)\mbox,
\end{equation}
with the same conditions as in Theorem \ref{thm:rigidaffine} and $\varepsilon:(\nC^d, 0) \rightarrow (\nC, 0)$.

Theorem \ref{thm:rigidaffine} says exactly then that we can kill $\varepsilon$. 
\end{myoss}

\begin{mydim}
Thanks to Remark \ref{oss:rigidd-1reducible}, we can suppose that $\ve{f}$ is of the form \eqref{eqn:rigidd-1reducible}.
We want to conjugate $f$ with a map $\widetilde{f}:(\nC^d, 0)\rightarrow (\nC^d, 0)$ of the form \eqref{eqn:rigidaffinepart} (with $\widetilde{\omega}$ instead of $\omega$).

We shall consider a local diffeomorphism of the form
$$
\ve{\Phi}(\vx)=\Big(\vxx, \vxy, \xz \big(1+\phi(\vx) \big)\Big)\mbox,
$$
where $\vx=(\vxx, \vxy, \xz)$, and $\phi:(\nC^d, 0) \rightarrow (\nC, 0)$.

Considering the conjugacy relation $\Phi \circ f = \widetilde{f} \circ \Phi$ for the last coordinate $\xz$, we get
\begin{align}
\xz \circ \ve{\Phi} \circ \ve{f}(\vx) &= \nu \vxx^\vi{l} \vxy^\vi{m} \xz \big(1+\varepsilon(\vx)\big)\big(1+\phi \circ \ve{f}(\vx)\big) + \omega(\vxx, \vxy) \big(1+\phi \circ \ve{f}(\vx)\big) \label{eqn:rigidaffineconjrel1}\\
\xz \circ \ve{\widetilde{f}} \circ \ve{\Phi}(\vx) &= \nu \vxx^\vi{l} \vxy^\vi{m} \xz \big(1+\phi(\vx)\big) + \widetilde{\omega}(\vxx, \vxy) \mbox. \nonumber
\end{align}
We want now to split \eqref{eqn:rigidaffineconjrel1} in two parts, one divisible by $\xz$, and the other that depends only on $(\vxx, \vxy)$.
Using the equivalence
\begin{equation}\label{eqn:affineintegralisdifference}
\int_{0}^1 \frac{d}{d \tau} \big(\phi \circ \ve{f}(\vxx, \vxy, \tau \xz)\big) = \phi\big(\ve{f}(\vxx, \vxy, \xz)\big)-\phi\big(\ve{f}(\vxx, \vxy, 0)\big)
\end{equation}
and by direct computation we get
\begin{align*}
\nu \vxx^\vi{l} \vxy^\vi{m} \xz &\left( \big(1+\varepsilon(\vx)\big)\big(1+\phi \circ \ve{f}(\vx)\big) + \omega(\vxx, \vxy)\int_0^1 \frac{\partial \phi}{\partial \xz}\big(\ve{f}(\vxx, \vxy, \tau \xz)\big) \big(1 + \zeta(\vxx, \vxy, \tau\xz)\big) d\tau\right)  \\
&+ \omega(\vxx, \vxy) \big(1+\phi \circ f (\vxx, \vxy, 0) \big)\mbox,
\end{align*}
where $\zeta: (\nC^d, 0) \rightarrow (\nC, 0)$ is given by
$$
\zeta(\vx):= \varepsilon(\vx) + \xz \frac{\partial \varepsilon}{\partial \xz}(\vx)\mbox.
$$
The conjugacy relation then gives two equations to solve (comparing the part divisible by $\xz$ and the one that does not depend on $\xz$), with respect to $\phi$ and $\widetilde{\omega}$:
\begin{align}
\varepsilon(\vx) + T\phi(\vx)&=\phi(\vx)\mbox, \label{eqn:affineeqntosolve1}\\
\omega(\vxx, \vxy) \big(1+\phi \circ f (\vxx, \vxy, 0) \big) &= \widetilde{\omega}(\vxx, \vxy) \mbox, \label{eqn:affineeqntosolve2}
\end{align}
where $\psi \mapsto T\psi$ is the functional given by
\begin{equation}\label{eqn:defT}
(T\psi)(\vx):= \big(1+\varepsilon(\vx)\big)\psi \circ \ve{f}(\vx) + \omega(\vxx, \vxy)\int_0^1 \frac{\partial \psi}{\partial \xz}\big(\ve{f}(\vxx, \vxy, \tau \xz)\big) \big(1 + \zeta(\vxx, \vxy, \tau\xz)\big) d\tau\mbox.
\end{equation}
Equation \eqref{eqn:affineeqntosolve1} has a solution given by
$$
\phi(\vx)=\sum_{n=0}^\infty T^{\comp n}\varepsilon(\vx)\mbox,
$$
we refer to the proof in \cite{favre:rigidgerms} in the $2$-dimensional case for convergence estimates, that rely on Cauchy's estimates.
Once that $\phi$ is defined as a holomorphic germ that satisfies \eqref{eqn:affineeqntosolve1}, we can use \eqref{eqn:affineeqntosolve2} to define $\widetilde{\omega}$, and we are done.
\end{mydim}

\begin{myoss}
Theorem \ref{thm:rigidaffine} tells us that, given a rigid germ $f: (\nC^d,0) \rightarrow (\nC^d, 0)$ of the form \eqref{eqn:rigidd-1reducible}, we can change coordinates holomorphically in order to have that the last coordinate of $\ve{f}$ is an affine function on $\xz$ (with coefficients that depend on the other coordinates $\vxx, \vxy$).
\end{myoss}

\begin{myoss}\label{oss:nonpolnormalform}
While studying rigid germs under the hypothesis of Theorem \ref{thm:rigidaffine}, following the argument used in the classification of $2$-dimensional contracting rigid germs (see \cite[pp. 494--498]{favre:rigidgerms}), we should consider change of coordinates of the form 
\begin{equation}\label{eqn:changecoordd-1nonpol}
\Phi(\vxx,\vxy,\xz) = \big(\vxx, \vxy, \xz + \phi(\vxx, \vxy)\big)
\end{equation}
(we are using the notations of Theorem \ref{thm:rigidaffine}).
In dimension $2$, one can obtain (holomorphically) that $\omega$ is a polynomial map in $(\vxx, \vxy)=\x_1$.
This is no longer true in general, not even formally, in higher dimensions.
Indeed, by computing the coefficients in the conjugacy relation, one can show that there can be infinitely many coefficients of $\omega$ that cannot be killed up to a change of coordinates of the form \eqref{eqn:changecoordd-1nonpol}.
It can be also shown that, in order to maintain the normal form as in \eqref{eqn:rigidaffinepart}, one can (basically) consider only change of coordinates such as \eqref{eqn:changecoordd-1nonpol}.
\end{myoss}

\begin{myoss}\label{oss:nonarchiaffine}
Theorem \ref{thm:rigidaffine} holds over any complete metrized field $\nK$ of characteristic $0$ (provided as always that all eigenvalues of $df_0$ belong to $\nK$, see Remark \ref{oss:nonarchijordan}).
Indeed, in the whole proof we never take roots of polynomials, so the argument works also for non-algebraically closed fields.

In the proof of \ref{thm:rigidaffine}, we define and estimate an operator $T$ given by \eqref{eqn:defT}.
To define $T$ we use integrals, so convergence could fail for the presence of (big) integers as denominators of the coefficients of the formal power series involved.
But thanks to \eqref{eqn:affineintegralisdifference}, we can write the integral appearing in \eqref{eqn:defT} as a difference of convergent formal power series.
Moreover, to prove convergence we use Cauchy's estimates, that are even stronger in the non-archimedean setting.
It follows that the argument works also for non-archimedean fields.
\end{myoss}

\mysect{Rigid Germs in Dimension $3$}
\begin{mytext}
With Table \ref{tab:normalformsd3} we summarize the normal forms obtained for a contracting rigid germ $f:(\nC^3, 0) \rightarrow (\nC^3, 0)$, with the assumption of injective internal action.
We set $q$ the number of irreducible components of $\critinf{f}$, $r$ the number of periodic components, $s$ the number of non-zero eigenvalues of $df_0$, $\eta \in \nN^*$ the order of (the matrix associated to) the periodic components of $\critinf{f}$.
We shall denote by $\mg{m}=\langle \X, \Y, \Z \rangle$ the maximal ideal of $\nC[[\X,\Y,\Z]]$. We shall also denote by $\lambda^1, \lambda^2, \lambda^3$ the eigenvalues of $df_0$ ordered as following:
\begin{equation*}
\abs{\lambda^1} \geq \abs{\lambda^2} \geq \abs{\lambda^3}\mbox.
\end{equation*}
\end{mytext}

\begin{table}
\caption{Contracting rigid germs for $d=3$.}
\label{tab:normalformsd3}
\begin{flushleft}
\begin{tabular}{|c|c|c|c|m{10.5cm}|}
\hline
$q$ & $r$ & $s$ & $\critinf{f}$ & Normal form\\
\hline
$0$ & $0$ & $3$ & $\emptyset$ & $\Big.\big(\lambda^1 \X, \lambda^2 \Y + \rho^1(\X), \lambda^3 \Z + \rho^2(\X, \Y)\big)\Big.$, the Poincar\'e-Dulac normal form.\\
\hline
$1$ & $0$ & $0$ & $\{\X=0\}$ & $\Big.\big(\beta \X^d, ?, ?\big)\Big.$, $d \geq 2$, $\beta \in \nC^*$.\\
\cline{3-5}
& & $1$ & $\{\Y=0\}$ & $\Big.\big(\lambda^1 \X, \Y^d, \nu \Y^m \Z + \omega(\X, \Y)\big)\Big.$, $d \geq 2$, $m \geq 1$, $\nu \in \nC^*$,  $\omega(\X, \Y) - \varepsilon \Y \in \mg{m}^2$ for a suitable $\varepsilon \in \{0,1\}$. \\
\cline{3-5}
& & $2$ & $\{\Z=0\}$ & $\Big.\big(\lambda^1 \X, \lambda^2 \Y + \rho \X^n, \Z^d\big)\Big.$, $\rho \in \{0, 1\}$ if $(\lambda^1)^n = \lambda^2$, $\rho = 0$ otherwise; $d \geq 2$. \\
\cline{2-5}
& $1$ & $1$ & $\{\X=0\}$ & $\Big.\big(\lambda^1 \X, ?, ?\big)\Big.$. \\
\cline{3-5}
& & $2$ & $\{\X=0\}$ & $\Big.\big(\lambda^1 \X, \lambda^2 \Y + \rho \X^n, \X^l \Z + \omega(\X, \Y)\big)\Big.$, $\rho \in \{0,1\}$ if $(\lambda^1)^n = \lambda^2$, $\rho = 0$ otherwise;\\
&  &  & $\{\Y=0\}$ & $\Big.\big(\lambda^1 \X, \lambda^2 \Y, \Y^l \Z + \omega(\X, \Y)\big)\Big.$;\\
&  &  &  & in both cases, $l \geq 1$, $\omega \in \mg{m}^2$.\\
\hline
$2$ & $0$ & $0$ & $\{\X\Y=0\}$ & $\Big.\big(\beta^1 \X^{d_1^1} \Y^{d_2^1}, \beta^2 \X^{d_1^2} \Y^{d_2^2}, \nu \X^l \Y^m \Z + \omega(\X, \Y)\big)\Big.$, $\beta^1, \beta^2 \in \nC^*$, $d_1^1 d_2^2 \neq d_1^2 d_2^1$, $d_1^2+d_2^2 \geq 2$, $\max\{d_1^1-1, d_2^1\} \geq 1$, $\nu \in \nC^*$, $l+m \geq 1$, $\omega \in \mg{m}^2$. \\
\cline{3-5}
& & $1$ & $\{\Y\Z=0\}$ & $\Big.\big(\lambda^1 \X, \beta^1 \Y^{d_1^1} \Z^{d_2^1} (1 + g \X^n), \beta^2 \Y^{d_1^2} \Z^{d_2^2}\big)\Big.$, $\beta^1, \beta^2 \in \nC^*$, $d_1^1 d_2^2 \neq d_1^2 d_2^1$, $\max\{d_1^1-1, d_2^1\} \geq 1$, $g \in \{0,1\}$ if $\big((\lambda^1)^n -d_1^1\big) \big((\lambda^1)^n - d_2^2\big)=d_1^2 d_2^1$, $g=0$ otherwise. \\
\cline{2-5}
& $1$ & $1$ & $\{\X\Y=0\}$ & $\Big.\big(\lambda^1 \X, \X^c \Y^d, \nu \X^l \Y^m \Z + \omega(\X, \Y)\big)\Big.$, $c+d \geq 2$, $l+m \geq 1$, $c+l \geq 1$, $d \geq 1$, $d+m \geq 2$, $\nu \in \nC^*$, $\omega(\X, \Y) - \varepsilon \Y \in \mg{m}^2$ for a suitable $\varepsilon \in \{0,1\}$. \\
\cline{3-5}
& & $2$ & $\{\X\Z=0\}$ & $\Big.\big(\lambda^1 \X, \lambda^2 \Y + \rho \X^n, \X^c \Z^d\big)\Big.$, $\rho \in \{0,1\}$ if $(\lambda^1)^n = \lambda^2$, $\rho = 0$ otherwise;\\
&  &  & $\{\Y\Z=0\}$ & $\Big.\big(\lambda^1 \X, \lambda^2 \Y, \Y^c \Z^d\big)\Big.$;\\
&  &  &  & in both cases, $c \geq 1$, $d \geq 2$.\\
\cline{2-5}
& $2$ & $2$ &  $\{\X\Y=0\}$ & \begin{minipage}{10cm}
$\eta = 1$: $\Big.\big(\lambda^1 \X, \lambda^2 \Y, \X^l \Y^m \Z + \omega(\X, \Y)\big)\Big.$;\\ 
$\eta = 2$: $\Big.\big(\alpha^1 \Y, \alpha^2 \X, \X^l \Y^m \Z + \omega(\X, \Y)\big)\Big.$, $\alpha^1 \alpha^2 = -\lambda^1 \lambda^2$;\\
in both cases, $l, m \geq 1$, $\omega \in \mg{m}^2$.
                              \end{minipage}\\
\hline
$3$ & $0$ & $0$ & $\{\X\Y\Z=0\}$ & $\Big.\big(\beta^1 \X^{d_1^1}\Y^{d_2^1}\Z^{d_3^1}, \beta^2 \X^{d_1^2}\Y^{d_2^2}\Z^{d_3^2}, \beta^3 \X^{d_1^3}\Y^{d_2^3}\Z^{d_3^3}\big)\Big.$, $\beta^1, \beta^2, \beta^3 \in \nC^*$, $\mi{D}:=(d_i^j)$ such that $\det \mi{D} \neq 0$, $d_1^j + d_2^j + d_3^j \geq 2$ for $j=1,2,3$. \\
\cline{2-5}
& $1$ & $1$ & $\{\X\Y\Z=0\}$ & $\Big.\big(\lambda^1 \X, \beta^1 \X^{c_1}\Y^{d_1^1}\Z^{d_2^1} (1+ g \X^n), \beta^2 \X^{c_2} \Y^{d_1^2} \Z^{d_2^2} \big)\Big.$, $\beta^1, \beta^2 \in \nC^*$, $d_1^1 d_2^2 \neq d_1^2 d_2^1$, $c_1 + c_2 \geq 1$, $d_1^j + d_2^j \geq 2$ for $j=1,2$, $g \in \{0,1\}$ if $\big((\lambda^1)^n -d_1^1\big) \big((\lambda^1)^n - d_2^2\big)=d_1^2 d_2^1$, $g = 0$ otherwise. \\
\cline{2-5}
& $2$ & $2$ & $\{\X\Y\Z=0\}$ & \begin{minipage}{10cm}
$\eta = 1$: $\Big.\big(\lambda^1 \X, \lambda^2 \Y, \X^{c_1} \Y^{c_2} \Z^d\big)\Big.$;\\ 
$\eta = 2$: $\Big.\big(\alpha^1 \Y, \alpha^2 \X, \X^{c_1} \Y^{c_2} \Z^d\big)\Big.$, $\alpha^1 \alpha^2 = -\lambda^1 \lambda^2$;\\
in both cases, $c_1, c_2 \geq 1$, $d \geq 2$.
                               \end{minipage}\\
\hline
\end{tabular}
 
\end{flushleft}

\end{table}

\begin{myoss}
By performing another change of coordinates of the form $(\X,\Y,\Z) \mapsto (\kappa^1 \X, \kappa^2 \Y, \kappa^3 \Z)$, with $\kappa^1, \kappa^2, \kappa^3 \in \nC^*$, we can say a little more on coefficients that arise in the normal forms.
\begin{itemize}
\item for $q=2$, $r=s=0$, we can put $2$ coefficients among $\beta^1, \beta^2, \nu$ equal to $1$ if the matrix
$$
\begin{mymatrix}{3}
d_1^1-1 & d_1^2 & l \\
d_2^1 & d_2^2-1 & m
\end{mymatrix}
$$
has rank $2$ (the ones associated to a $2 \times 2$ invertible submatrix), otherwise we can just put one of them equal to $1$ (for example $\nu=1$).
\item for $q=2$, $r=s=1$, we can put $\nu = 1$ if
$$
\det
\begin{mymatrix}{2}
c & l \\
d-1 & m 
\end{mymatrix}
\neq 0\mbox.
$$
\item for $q=2$, $r=0$, $s=1$, or $q=3$, $r=s=0,1$, if we put $D=(d_i^j)$, then we can put $\beta^j=1$ for as many $j$ as the rank of $D- \on{Id}$.
\end{itemize}
\end{myoss}

In this classification, two cases are not completely understood: $q=1$ and $r=s=0,1$, i.e., when $p+s=1$.
If we consider the action of $f$ on $\critinf{f}=\{\X=0\}$, we can have two behaviors: either $f(\{\X=0\})=0$, or $f(\{\X=0\})$ is a (not necessarily smooth) curve in $\{\X=0\}$.
The following example will show that this second case can happen for every irreducible curve in $\{\X=0\}$.

\begin{myes}\label{es:anycurve}
Let $\Psi:(\nC, 0) \rightarrow (\nC^2, 0)$ be the parametrization of a curve $\mc{C}$, of the form
$$
\Psi(\T)=\big(\T^m, \psi(\T)\big)\mbox,
$$
where $\psi:(\nC,0) \rightarrow (\nC,0)$ is a holomorphic map with multiplicity $m(\psi) \geq m$ at $0$.

Consider the map $f:(\nC^3, 0) \rightarrow (\nC^3, 0)$ given by
$$
(\X, \Y, \Z) \mapsto \big(\lambda \X^a, \X\Y + \Z^m, \X\Z + \X\Y \xi(\Z) + \psi(\Z)\big)\mbox,
$$
where $a \geq 1$, $\lambda \in \nC^*$ (and $\abs{\lambda} < 1$ if $a = 1$ to have a contracting germ), and $\xi:(\nC,0) \rightarrow \nC$ is given by
$$
\xi(\Z):=\frac{\psi^\prime(\Z)}{m \Z^{m-1}}\mbox.
$$
Computing the Jacobian, we get
$$
\det df = \lambda a \X^{a+1} \big(1+\Y \xi^\prime(\Z)\big)\mbox,
$$
and hence $f$ is a contracting rigid germ such that $f(\{\X=0\})=\mc{C}$.
\end{myes}

\begin{mytext}
Example \ref{es:anycurve} shows how, to study the classification of the missing cases, we have to take care of the geometry of the images of $\critinf{f}$, and maybe make some additional assumptions to get some classification results.

With the next example, we shall show another phenomenon that can appear.
\end{mytext}

\begin{myes}\label{es:manyimages}
Consider the map $f:(\nC^3, 0) \rightarrow (\nC^3, 0)$ given by
$$
(\X, \Y, \Z) \mapsto \big(\lambda \X^a, \X(1+\Y^2), \X\Y\Z^2\big)\mbox,
$$
where $a \geq 1$ and $\lambda \in \nC^*$ (and $\abs{\lambda} < 1$ if $a = 1$ to have a contracting germ).

Then $\critinf{f}=\{\X\Y\Z=0\}$, while
\begin{align*}
f(0,\Y,\Z)&=(0,0,0)\mbox,\\
f(\X, 0, \Z)&=(\lambda \X^a, \X, 0)\mbox,\\
f(\X, \Y, 0)&=\big(\lambda \X^a, \X(1+\Y^2), 0\big)\mbox, 
\end{align*}
hence $f(\critinf{f})\subseteq \{\Z=0\} \subset \critinf{f}$, and $f$ is rigid.

But by direct computation we get that $f^{\comp n}(\{\Y=0\})=:\mc{C}_n$ form a sequence of distinct curves in $\{\Z=0\} \cong (\nC^2, 0)$.

The geometry of $\bigcup_n \mc{C}_n$, or rather of $\Delta \setminus \bigcup_n \mc{C}_n$ where $\Delta$ is a small polydisc centered in $0$, should be taken into account to find a classification up to holomorphic (or even formal) change of coordinates.  
\end{myes}

\bibliographystyle{alpha}
\bibliography{biblio}

\vspace{0.5cm}

{\small \noindent
Matteo Ruggiero,\\
Fondation Math\'ematique Jacques Hadamard,
D\'epartement de Math\'ematiques,
UMR 8628 Universit\'e Paris-Sud 11-CNRS,
B\^atiment 425,
Facult\'e des Sciences d'Orsay,
Universit\'e Paris-Sud 11,
F-91405 Orsay Cedex, France.\\
Centre de Math\'ematiques Laurent Schwartz, 
\'Ecole Polytechnique, 
91128 Palaiseau Cedex, France.\\
Tel: (+33) (0)1 69 33 49 18. E-mail: ruggiero@math.polytechnique.fr \\
}

\end{document}